\documentclass[11pt]{article}
\usepackage{amssymb,amssymb,amscd}
\begin{document}
\newtheorem{thm}{Theorem}
\newtheorem{lem}[thm]{Lemma}
\newtheorem{cor}[thm]{Corollary}
\newtheorem{prop}[thm]{Proposition}
\def\endpf{\hfill$\bullet$\medskip}
\newtheorem{rem}[thm]{Remark}
\newtheorem{exam}[thm]{Example}
\newtheorem{defn}[thm]{Definition}
\newcounter{jolist}
\newenvironment{joliste}{\begin{list}{{\rm (\roman{jolist})}}{%
    \usecounter{jolist}\setlength{\labelwidth}{10mm}
    \setlength{\leftmargin}{10mm}\setlength{\itemsep}{-4pt}
   \setlength{\topsep}{0pt}}}{\end{list}}
\newcommand{\al}{\alpha}
\newcommand{\avee}{\alpha^\vee}
\newcommand{\bvee}{\beta^\vee}
\newcommand{\wt}{\widetilde}
\newcommand{\delp}{\partial^+}
\newcommand{\delm}{\partial^-}
\newcommand{\lam}{{\lambda}}
\newcommand{\Lam}{{\Lambda}}
\newcommand{\co}{\mathcal{O}}
\newcommand{\ck}{\mathcal{K}}
\newcommand{\cF}{\mathcal{F}}
\newcommand{\cl}{\mathcal{L}}
\newcommand{\ci}{\mathcal{I}}
\newcommand{\cf}{\mathcal{F}}
\newcommand{\cp}{mathcal{P}}
\newcommand{\charc}{\mathop{\rm Char}}
\newcommand{\Ker}{\mathop{\rm Ker}}
\newcommand{\Lie}{\mathop{\rm Lie}}
\newcommand{\Pic}{\rm Pic}
\newcommand{\Ch}{\rm Ch}
\newcommand{\res}{\mathop{\rm res}}
\newcommand{\proj}{\mathop{\rm proj}}
\newcommand{\grad}{\mathop{\rm grad}}
\newcommand{\bp}{{\mathbb P}}
\newcommand{\bn}{{\mathbb N}}
\newcommand{\bz}{{\mathbb Z}}
\newcommand{\bq}{{\mathbb Q}}
\newcommand{\bb}{{\mathbb B}}
\newcommand{\bd}{{\mathbb D}}
\newcommand{\usigma}{{\underline{\sigma}}}
\newcommand{\udelta}{{\underline{\delta}}}
\newcommand{\ukappa}{{\underline{\kappa}}}
\newcommand{\ua}{\underline{a}}
\newcommand{\ub}{\underline{b}}
\newcommand{\uc}{\underline{c}}
\newcommand{\ui}{{\underline{i}}}
\newcommand{\uj}{{\underline{j}}}
\newcommand{\uka}{{\underline{k}}}
\newcommand{\ut}{{\underline{\tau}}}
\newcommand{\uk}{{\underline{\kappa}}}
\newcommand{\upi}{{\underline{\pi}}}
\newcommand{\ueta}{{\underline{\eta}}}
\newcommand{\umu}{{\underline{\mu}}}
\newcommand{\unu}{{\underline{\nu}}}
\newcommand{\us}{{\underline{\sigma}}}
\newcommand{\tbar}{{\overline{\tau}}}
\newcommand{\kbar}{{\overline{\kappa}}}
\newcommand{\sbar}{{\overline{\sigma}}}
\newcommand{\ts}{\tilde{\sigma}}
\newcommand{\cR}{{\hat{R}}}
\newcommand{\km}{{\mathfrak{g}}}
\newcommand{\lib}{{\mathfrak{b}}}
\newcommand{\liu}{{\mathfrak{n}}}
\newcommand{\libmin}{{\mathfrak{b}^-}}
\newcommand{\liumin}{{\mathfrak{n}^-}}
\title{Richardson Varieties and\\ Equivariant $K$-Theory}
\author{V. Lakshmibai\footnote{V.L. supported by NSF Grant-DMS-9971295.} and
P. Littelmann\footnote{P.L. supported by TMR-Grant ERB FMRX-CT97-0100.}}
\maketitle
\begin{abstract}
We generalize Standard Monomial Theory (SMT) to intersections of
Schubert varieties and opposite Schubert varieties; such varieties
are called Richardson varieties. The aim of this article is to get
closer to a geometric interpretation of the standard monomial
theory as constructed in \cite{Li$_2$}. In fact, the construction
given here is very close to the ideas in \cite{$G/P$-V}. Our
methods show that in order to develop a SMT for a certain class of
subvarieties in $G/B$ (which includes $G/B$), it suffices to have
the following  three ingredients, a basis for $H^0(G/B,\cl_\lam)$,
compatibility of such a basis with the varieties in the class,
certain quadratic relations in the monomials in the basis
elements. An important tool (as in \cite{$G/P$-V}) will be the
construction of nice filtrations of the vanishing ideal of the
boundary of the varieties above. This provides a direct connection
to the equivariant $K$-theory (products of classes of structure
sheaves with classes of line bundles), where the combinatorially
defined notion of standardness gets a geometric interpretation.
\end{abstract}

\section*{Introduction}\label{Introduction}

In many respects, the simplest case of the construction of a
Standard Monomial basis is the case of the Grassmannian $X$ of
linear subspaces of dimension $r$ in $k^n$, embedded into the
projective space $\bp(\wedge^r k^n)$ via the Pl\"ucker embedding;
let $L$ be the corresponding very ample line bundle on $X$. Then
the homogeneous coordinate ring of the embedded variety coincides
with the ring $R=\bigoplus_{m=0}^{\infty} H^0(X,L^{\otimes m})$,
and this ring admits a nice basis, defined as follows. Let
$\{e_1,\ldots,e_n\}$ be the standard basis of $k^n$, then the
wedge products $e_\ui=e_{i_1} \wedge \cdots \wedge e_{i_r}$, where
$\ui$ varies over the set $I$ of all $r$-tuples such that  $1\le
i_1<\ldots<i_r\le n$, form a basis of $\wedge^r k^n$.

Let $\bb=\{p_\ui\mid\ui\in I\}$ be the dual basis. These elements
form naturally a generating set of the homogeneous coordinate ring
of the embedded variety, and it is a natural question to ask for a
description of the ring in terms of relations, and, similarly, for
the subvarieties like Schubert varieties, opposite Schubert
varieties, unions and intersections of these. To formulate this
more precisely, first observe that we have a natural partial order
on $I$: $\ui\ge\uj$ if $i_1\ge j_1,\ldots, i_{r}\ge j_{r}$.

A product of elements $p_\ui p_\uj \ldots p_\uka$ is called {\it
standard} if the indices are linearly ordered with respect to this
partial order, i.e., $\ui\ge\uj\ge\ldots\ge\uka$. These monomials
deserve this ``special name'' because non-standard products can be
expressed as linear combinations of standard products. This
follows from the well-known Pl\"ucker relations: if two elements
$p_\ui$ and $p_{\ui'}$ are not comparable in the partial order,
then the following relation, called {\it quadratic straightening
relation}, holds in the homogeneous coordinate ring: $$ p_\ui
p_{\ui'}= \sum_{\uj\ge\ui,\ui'\ge\uka}a_{\uj,\uka}p_\uj p_{\uka},
$$ where the coefficients are elements of the ground field $k$ and
the products $p_\uj p_{\uka}$ occurring on the right side of the
equation are standard. Hodge \cite{H} has already observed that
these relations give a full presentation of the ring $R$. In
particular, the standard monomials form a basis of the ring $R$.
Further, Hodge has also noticed that this basis is compatible with
Schubert varieties, i.e., for $\ui\in I$ let $X_\ui$ be the
Schubert variety consisting of subspaces $U\in X$ such that for
all $1\le s\le r$, $\dim(U\cap\hbox{\rm span}(e_1,\ldots,e_s))$ is
greater or equal to the number of $j$ such that $i_j\le s$. We say
that a standard monomial $p_\uj \ldots p_\uka$ is standard on
$X_\ui$ if $\ui\ge \uj$. Then Hodge has already shown that the
restriction of the standard monomials, not standard on $X_\ui$,
vanish identically, and the restrictions of the standard
monomials, standard on $X_\ui$, form a basis of
$R_\ui=\bigoplus_{m=0}^{\infty} H^0(X_\ui,L^{\otimes m})$. Such a
description has some immediate geometric applications, for example
the Schubert variety is projectively normal in the embedding, one
has a presentation of the homogeneous ideal of the Schubert
varieties etc.

The purpose of {\it Standard Monomial Theory} (SMT) is to generalize this kind of description
to all embeddings $G/P\hookrightarrow \bp(H^0(G/P,\cl_\lam))$, where $G$ is a semi-simple
algebraic group over an algebraically closed field, $P$ is a parabolic subgroup,
and $\cl_\lam$ is an ample line bundle on $G/P$ corresponding to a dominant weight
of $G$.

This paper has two aims: The first aim is to generalize SMT to intersections of Schubert
varieties and opposite Schubert varieties; such varieties will be called Richardson
varieties, since Richardson looked at these first (cf. \cite{Ri}). The second aim of this paper
is to make SMT easier generalizable to other classes of varieties.

To be more precise, fix a Borel subgroup $B\subset G$ and a
maximal torus $T$, and let $B^-$ denote the Borel subgroup of $G$
opposite to $B$ (it is the unique Borel subgroup of $G$ with the
property $B\cap B^-=T$). Let $P\supset B$ be a parabolic subgroup,
let $W_P$ be its Weyl group, we identify the group canonically
with a subgroup of the Weyl group $W$ of $G$. For $\tau \in
W/W_P$, let $e_\tau\in G/P$ be the corresponding coset and let
$X_\tau$ denote the {\it Schubert variety}, the Zariski closure of
the $B$-orbit $B e_ \tau$ in $G/P$. Similarly, let $X^\tau$ denote
the {\it opposite Schubert variety}, the Zariski closure of the
$B^-$-orbit $B^-e_ \tau$ in $G / P$. The {\it Richardson variety}
$X_\tau^\kappa$ is then the intersection $X_\tau\cap X^\kappa$.
For $P=B$, such double coset intersections $B\tau B\cap B^-\kappa
B$ first appear in \cite{Bi}, \cite{Deo}, \cite{K-Lu$_1$},
\cite{K-Lu$_2$}, \cite{Ri}, \cite{Ri-Ro-S}; in fact, it is shown
in  \cite{Ri} that $B\tau B\cap B^-\kappa B$ is dense in
${\overline{B\tau B}}\cap {\overline{B^-\kappa B}}$. Recently,
such (as well as similar) double coset intersections have appeared
in the context of total positivity (cf. \cite{Be-F-Z},
\cite{Be-Z}, \cite{Fo-Z}, \cite{Lu$_1$}, \cite{Rie$_1$},
\cite{Rie$_2$}, \cite{Rie$_3$}, \cite{S-S-V}). Richardson
varieties also appear in the context of K-theory of flag varieties
as explained below.

The starting point of our approach will be as in the example
above: We start with

{\it a)}  a basis $\bb(\lam)$ of $H^0(G/P,\cl_\lam)$ which will be
indexed by a partially ordered set $B(\lam)$

{\it b)}  $\bb(\lam)$ is compatible with Schubert and opposite
Schubert varieties $Z$ (i.e., the set $\{b\vert_Z\mid b\in
\bb(\lam), b\vert_Z\not\equiv 0\}$ is linearly independent)

 {\it c)} $\bb(\lam)$ satisfies certain quadratic relations similar to the {\it quadratic
straightening relations} described above.

The main step in the construction of SMT is the following: let
$X_\tau^\kappa$ be a Richardson variety. By the boundary $\delp
X_\tau^\kappa$ we mean the subvariety of $X_\tau^\kappa$
consisting of the union of all Richardson varieties
$X_{\tau'}^{\kappa'}$ such that $\kappa\le\kappa'\le\tau'<\tau$.
Let $\delp\ci_\tau^\kappa$ be the ideal sheaf of this subvariety.
Using the global basis $\bb(\lam)\subset H^0(G/P,\cl_\lam)$ and
the straightening relations, we show that the twisted sheaf
$\delp\ci_\tau^\kappa\otimes\cl_\lam$ admits a filtration as a
coherent $\co_X$-sheaf such that the associated graded
$\co_X$-sheaf is a direct sum of structure sheaves of Richardson
varieties, and for the coefficients we have a combinatorial
formula involving the combinatorics of the indexing set $B(\lam)$.
Once we have such a {\it Pieri--Chevalley type formula}, we can
proceed in a straight forward way (using induction on the
dimension) and establish a standard monomial theory compatible
with all unions of Richardson varieties.

Summarizing, we conclude the introduction with a short description
of the development of SMT and its present status. It was developed
by Lakshmibai, Musili, and Seshadri in a series of papers,
culminating in \cite{$G/P$-V} where it is established for all
$G/P\hookrightarrow \bp(H^0(G/P,\cl_\lambda))$, $G$ a classical
group. Further results concerning certain exceptional and
Kac--Moody groups led to conjectural formulations of a general
SMT, see \cite{LS2}. The combinatorial conjectures have been
proved by the second author in \cite{Li$_1$}, where he introduced
a new combinatorial tool, the path model. Another development was
the introduction of the dual of Lusztig's Frobenius map for
quantum groups at roots of unity (see \cite{Li$_2$}, see also
\cite{KL1,KL2}). This Frobenius map turns out to provide exactly
the right tool to give a construction of SMT in the most general
setting: $G/P\hookrightarrow \bp(H^0(G/P,\cl_\lam))$, where $G$ is
an arbitrary symmetrizable Kac--Moody group and $\lam$ is an
integral dominant weight such that $\cl_\lam$ is ample on $G/P$.

To be able to generalize SMT to other types of varieties (compactifications of symmetric spaces or, more
generally, spherical varieties for example), it would be very helpful to get a better algebraic geometric
description or characterization of the SMT. A first step in this direction is the observation of the
connection between SMT and equivariant $K$-theory as described in~\cite{LiSe}, and this
article can be described as an effort to reduce the input from quantum group methods (which are not
available in the examples above) in \cite{Li$_2$} as much as possible and to use the Pieri--Chevalley type
formula as part of the construction instead of proving it as a consequence of SMT (as in~\cite{LiSe}). In
fact, once one has a good candidate for SMT, the inductive procedure and a Pieri--Chevalley type formula reduces
the proof to the counting of global sections while the proof of the vanishing theorem for higher cohomology has
reminiscences of the corresponding proofs using characteristic $p>0$ methods (see Theorem~\ref{vanishandbasis}).

Another step in this direction has been proved recently by M.
Brion and the first author (cf. \cite{brionlak}). They provide a
geometric construction of a (SMT) basis for $H^0(G/P,\cl_\lam)$,
compatible with unions of Richardson varieties. The main tool
introduced in \cite{brionlak} is a flat $T$-invariant degeneration
of the diagonal in $X\times X$ ($X$ being $G/P$) to the union of
all products $X_\tau\times X^\tau$, where the $X_\tau$ are the
Schubert varieties and the $X^\tau$ are the corresponding opposite
Schubert varieties; we also obtain this degeneration as a
consequence of the standard monomial basis, see
section~\ref{SMTequivariantKtheory}. It is interesting to note
that while in  \cite{brionlak}, the above degeneration is the
starting point towards developing a SMT for Richardson varieties,
in this paper such a degeneration is obtained as a consequence of
SMT ! In \cite{brionlak}, for any Richardson variety
$X_\tau^\kappa$, the authors first construct a basis of the space
of sections of $\cl_\lam$ on $X_\tau^\kappa$ vanishing on both of
the boundaries $\partial^{\pm}X_\tau^\kappa$ (here,
$\partial^{-}X_\tau^\kappa$ is the subvariety of $X_\tau^\kappa$
consisting of the union of all Richardson varieties
$X_{\tau'}^{\kappa'}$ such that $\kappa<\kappa'\le\tau'\le\tau$).
This basis is then lifted to global sections of $\cl_\lam$ that
vanish on all Richardson varieties not containing $X_\tau^\kappa$,
which in turn give raise to a basis for $H^0(G/P,\cl_\lam)$; this
basis for $H^0(G/P,\cl_\lam)$ is then shown to have compatibility
with unions of Richardson varieties. In this article, we show the
compatibility of the basis of \cite{Li$_2$} with unions of
Richardson varieties (cf. Theorem \ref{vanishandbasis})), and the
vanishing of a typical basis vector on
$\partial^{\pm}X_\tau^\kappa$ for suitable $\kappa$ and $\tau$.
Thus the basis as constructed in \cite{Li$_2$} may be considered
as a specially nice basis of the type constructed in
\cite{brionlak}. In spite of this relationship between the bases
of \cite{brionlak} and \cite{Li$_2$}, there are still some missing
links; to be very precise, let $\pi=(\ut,\ua)$ be an L-S path of
shape $\lam$, where $\ut=(\tau_0,\ldots,\tau_r)$ is a strictly
decreasing sequence in $W/W_P$. The first and last element
describe the smallest Richardson variety $X_{\tau_0}^{\tau_r}$
such that the restriction of the path vector $p_\pi$ does not
vanish identically. The other terms come up naturally in the
quantum group construction and have a suggestive algebraic
geometric interpretation in \cite{Li$_2$}, but they do not have
yet an interpretation in the construction in \cite{brionlak}.

The paper is organized as follows: In the first three sections, we
recall the combinatorics of the path model and the construction of
the path vector basis: the indexing set $I$ in the Grassmannian
case will be replaced by the set of L-S paths $B(\lam)$ of shape
$\lam$ (section~\ref{pathmodel}), the basis $\{p_\ui, \ui \in I\}$
will be replaced by the basis $\bb(\lam)$ consisting of path
vectors $\{p_\pi, \pi\in B(\lam)\}$ (section~\ref{pathvector}),
and the Pl\"ucker relations will be replaced by the quadratic
straightening relations (section~\ref{smtsection}). In
section~\ref{gradedringsandmodules} we recall the connection
between coherent sheaves and graded finitely generated modules
over the homogeneous coordinate ring $k[G/P]$ given by the
embedding. Using the properties {\it (a)--(c)}, we show in
section~\ref{ideal} that the ideal $I_\tau^\kappa$ has a basis
given by standard monomials. This is used in
section~\ref{idealsheaves} to prove the Pieri--Chevalley type
formula for the ideal sheaf $\delp\ci_\tau^\kappa$, and, as a
consequence, we obtain in section~\ref{smtandvan} the construction
of SMT and the vanishing of higher cohomology groups. In
sections~\ref{nonregular} and \ref{SMTnonreg} we discuss the case
where the bundle $\cl_\lam$ is only base point free but not
necessarily ample. In the last
section~\ref{SMTequivariantKtheory}, the relationship between
K-theory and SMT is brought out, especially, the computation of
the coefficients of the classes of structure sheaves appearing in
the product of the class of the structure sheaf of a Schubert
variety with the class of a line bundle.

The first author is grateful to the Universit\"at  Wuppertal for
the hospitality extended during her visit in May 2001; it was
there that this work originated.

\section{Notation}\label{Notation}
In this section we fix some standard notation which will be used
throughout the paper. The ground field $k$ is supposed to be
algebraically closed of arbitrary characteristic. The group $G$ is
a semisimple, simply connected algebraic group defined over $k$.
We fix a Borel subgroup $B\subset G$ and a maximal torus $T\subset
B$. We denote by $B^-$ the opposite Borel subgroup, i.e., $B^-$ is
the unique Borel subgroup of $G$ such that $B\cap B^-=T$. The
unipotent radical of $B$ is denoted by $U$, and the unipotent
radical of $B^-$ is denoted by $U^-$.

We use the same notation (but with gothic letters) for the corresponding Lie algebras, i.e.,
$\km$ is the Lie algebra of $G$, $\lib$ the Lie algebra of $B$, $\libmin$ is the Lie
algebra of $\Lie B^-$ and $\liu^\pm$ denotes the Lie algebra of $U^\pm$.
The corresponding enveloping algebras are denoted by $U(\km),U(\lib^\pm)$ and $U(\liu^\pm)$.

Let $Q\supset B$ be a parabolic subgroup of $G$ containing $B$.
The projective variety $G/Q$ admits a finite number of $B$-orbits
(as well as a finite number of $B^-$-orbits), they are indexed by
the $T$-fixed points in $G/Q$. Let $W$ be the Weyl group of $G$
and let $W_Q$ be the Weyl group of $Q$. The group $W_Q$ can be
canonically identified with the subgroup of $W$ generated by the
simple reflections $s_\al$ such that $-\al$ is a root of the root
system of $Q$. For $\tau\in W/W_Q$ let $e_\tau\in G/Q$ be the
corresponding $T$-fixed point. The closure of the $B$-orbit:
$X_\tau=\overline{B e_\tau}\subset G/Q$ is the {\it Schubert
variety} corresponding to $\tau$, and the closure of the
$B^-$-orbit: $X^\tau=\overline{B^- e_\tau}\subset G/Q$, the {\it
opposite Schubert variety} corresponding to $\tau$.

Let $\tau,\sigma\in W/W_Q$. The {\it Richardson variety}
corresponding to the pair $(\tau,\sigma)$ is the (set theoretic)
intersection $X_\tau^\sigma = X_\tau\cap X^\sigma$ (with the
induced reduced structure). Note taht if $\overline{w}_0\in W/W_Q$
is the class of the longest element $w_0$ in $W$, then
$X_{\overline{w}_0}^\sigma=X^\sigma$. Similarly, one has
$X_\tau^{id}=X_\tau$.

If $\uk=(\kappa_1,\ldots,\kappa_s)$ and $\us=(\sigma_1,\ldots,\sigma_s)$ are sequences
of elements in $W/W_Q$, then we denote by $X_{\uk}^{\us}$ the union $X_{\kappa_1}^{\sigma_1}
\cup\ldots\cup X_{\kappa_s}^{\sigma_s}$ of Richardson varieties with the induced reduced structure.

Let $\Lambda$ be the character group of $T$. We denote by
$\Lambda^+$ the dominant weights and by $\Lambda^{++}$ the regular
dominant weights. For the parabolic subgroup $Q$, let $\Lambda_Q$
be the subgroup of weights which can be (trivially) extended to
characters of $Q$, let $\Lambda^{+}_Q\subset \Lambda_Q$ be the
dominant ones and denote by $\Lambda^{++}_Q$ be the $Q$-regular
dominant weights, i.e., these dominant weights can not be extended
to characters of a parabolic subgroup containing $Q$ properly.

For $\lam\in \Lambda$, let $\cl_\lam=G\times_B k_{-\lam}$ be the
line bundle on $G/B$ associated to the $B$-character $(-\lam)$. A
geometric way to characterize $\Lambda_Q\subset \Lambda$ is to say
that $\lam$ is an element of $\Lambda_Q$ if and only if $\cl_\lam$
``goes down'' to a line bundle $\cl_\lam=G\times_Q k_{-\lam}$ on
$G/Q$,  and $\lam\in \Lambda_Q^{+}$ if only if $\cl_\lam$ is base
point free on $G/Q$, and $\lam\in \Lambda_Q^{++}$ if only if
$\cl_\lam$ is an ample line bundle on $G/Q$.

\section{The path model and some partial orders}\label{pathmodel}
An important combinatorial tool will be the L-S paths of shape
$\lam$, we recall quickly the most important features of the path
model. Let $\lam\in \Lam^+$ be a dominant weight. Fix a {\it total
order ``$\succeq$''} on $W/W_\lam$ refining the {\it Bruhat order
``$\ge$''}. Let $\ut=(\tau_0\succ\ldots\succ\tau_r)$ be a strictly
decreasing sequence of elements of $W/W_\lam$ and let
$\ua=(0<a_1<\ldots<a_r<1)$ be strictly increasing sequence of
rational numbers. The pair $\pi=(\ut,\ua)$ is called  a {\it
convex subset of shape $\lam$} of the orbit $W.\lam$. To motivate
the name, set $a_0=0$ and $a_{r+1}=1$ and set $x_i=a_i-a_{i-1}$
for $1\le i\le r+1$. Then $x_i$ is called the {\it weight of
$\tau_{i-1}$}, and the sum $\sum_{i=0}^r x_{i+1}\tau_{i}(\lam)$ is
a convex linear combination, called the {\it weight of $\pi$}, and
is denoted $\pi(1)$.

Such a {\it convex subset} $\pi$ is called a L-S path of shape $\lam$ if the sequence
$\ut=(\tau_0,\ldots,\tau_r)$ is strictly decreasing in the Bruhat order (on $W/W_\lam$), and
if the pair satisfies the following integrality condition. For all $i=1,\ldots,r$:
\begin{itemize}
\item set ${s_i}=l(\tau_{i-1})-l(\tau_i)$. There exists a sequence $\beta_1,\ldots,\beta_{s_i}$
of positive roots joining $\tau_{i-1}$ and $\tau_i$ by the corresponding reflections, i.e.,
\[
\tau_{i-1} >s_{\beta_1}\tau_{i-1}>s_{\beta_2}s_{\beta_1}\tau_{i-1}>\ldots
>s_{\beta_{s_i}}\cdots s_{\beta_1}\tau_{i-1}=\tau_i,
\]
and $a_i\langle\tau_{i-1}(\lam),s_{\beta_1}\cdots s_{\beta_{j-1}}(\beta^\vee_j)\rangle\in\bz$ for all
$j=1,\ldots,s_i$. Here  $\beta^\vee$ denotes the coroot of $\beta$, and $\langle\mu,\beta^\vee\rangle$
stands for the evaluation of $\mu\in\Lambda$ on the coroot.
\end{itemize}

For more details on the combinatorics of L-S paths we refer the
reader to \cite{Li$_1$}. Let $B(\lam)$ be the set of L-S paths of
shape $\lam$. The character of the Weyl module $V(\lam)$ of
highest weight $\lam$ can be calculated using the L-S paths:

\begin{thm}[\cite{Li$_1$}]\label{pathcharacter} $\charc V(\lam)=\sum_{\pi\in B(\lam)} e^{\pi(1)}$.
\end{thm}

Let $\pi=(\ut,\ua)$ be an L-S path of shape $\lam$, where
$\ut=(\tau_0,\ldots,\tau_r)$. We call $i(\pi)=\tau_0$, the {\it
initial direction} and $e(\pi)=\tau_r$, the {\it final direction}
of the path.

\begin{defn}\rm
Let $B(\lam)$ be the set of L--S paths of shape $\lam$. We say
that $\pi\in B(\lam)$ is {\it standard on a Richardson variety}
$X_\tau^\kappa$ if $\tau\ge i(\pi)$ and $e(\pi)\ge \kappa$, and we
denote by $B_\tau^\kappa(\lam)$ the set of all L--S paths of shape
$\lam$, standard on $X_\tau^\kappa$. If $\tau$ (resp. $\kappa$) is
the class of the longest word in $W$ (resp. $id$), then $\tau$
(resp. $\kappa$) will be omitted and we will write just
 $B^\kappa(\lam)$ ( resp. $B_\tau(\lam)$).
\end{defn}

A sequence $\upi=(\pi_1,\ldots,\pi_m)$ of
L-S paths of shape $\lam$ is called {\it standard} if
\begin{eqnarray}\label{standardsequence}
e(\pi_1)\ge i(\pi_2)\ge\ldots \ge e(\pi_{m-1})\ge i(\pi_m).
\end{eqnarray}
The notion of initial and final direction generalizes to standard sequences of length $m$
as follows: We set $i(\upi)=i(\pi_1)$ and $e(\upi)=e(\pi_m)$. The notion of a standard sequence
on a Richardson variety generalizes in the obvious way.

Let $\upi(1)=\pi_1(1)+\ldots+\pi_m(1)$ be the {\it weight of such a standard sequence}. The
theorem above generalizes to standard sequences:

\begin{thm}[\cite{Li$_1$}]\label{standardpathseq}
$\charc V(m\lam)=\sum e^{\upi(1)}$, where the sum runs over all standard sequences of
length $m$ of L-S paths of shape $\lam$.
\end{thm}

We need several orderings on the set of convex subsets of $W.\lam$.
Let $\pi=(\ut,\ua)$ and $\eta=(\uk,\ub)$ (where $\uk=(\kappa_0,\ldots,\kappa_s)$ and $\ub=(b_1,\ldots,b_s)$)
be two convex subsets of $W.\lam$. Induced by the Bruhat order on $W/W_\lam$, we have two types of partial
orders on such convex subsets.
\begin{itemize}
\item We say $\pi\ge\eta$ if $\pi$ is greater than $\eta$ in the weighted lexicographic sense:
$\tau_0 >\kappa_0$ or $\tau_0=\kappa_0$ and $a_1>b_1$ or $\tau_0=\kappa_0, a_1=b_1$ and
$\tau_1>\kappa_1$, etc.
\item We say $\pi\ge^r\eta$ if $\pi$ is greater than $\eta$ in the reverse weighted lexicographic sense:
$\tau_r >\kappa_s$ or $\tau_r=\kappa_s$ and $1-a_r>1-b_s$
or $\tau_r=\kappa_s, 1-a_r=1-b_s$ and $\tau_{r-1}>\kappa_{s-1}$, etc.
\end{itemize}

For a total order $\succeq$ on $W/W_\lam$ we define in the same
way total orders $\succeq$ and $\succeq^r$ on the set of convex
subset of elements of $W/W_\lam$.

We extend these (partial) orders in the obvious way to sequences. Let
$\upi=(\pi_1,\ldots,\pi_t)$ and $\ueta=(\eta_1,\ldots,\eta_t)$ be two
sequences of weighted subsets.
\begin{itemize}
\item We say $\upi\ge \ueta$ if $\pi_1>\eta_1$ or $\pi_1=\eta_1$ and $\pi_2>\eta_2$ etc.
\item We say $\upi\ge^r \ueta$ if $\pi_t>\eta_t$ or $\pi_t=\eta_t$ and $\pi_{t-1}>\eta_{t-1}$ etc.
\end{itemize}

The total orders $\succeq$ and $\succeq^r$ are extended to sequences in the same way.

\section{The path vector basis}\label{pathvector}
Suppose $\lam\in \Lambda_Q^{+}$, then the line bundle $\cl_\lam$
is base point free. Let $V(\lam)$ be the Weyl module of highest
weight $\lam$. The parabolic subgroup $Q$ stabilizes the line
through the highest weight vector $v_\lam$ and we have a
corresponding map $$ G/Q \rightarrow \bp(V(\lam)),\quad gQ\mapsto
[gv_\lam]. $$ Further, $H^0(G/Q,\cl_\lam)=V(\lam)^*$ is the dual
space of $V(\lam)$. Associated to the combinatorial path model
given by the L-S paths $B(\lam)$ of shape $\lam$, we have the
basis $$ \bb(\lam)=\{p_\pi\mid\pi\in B(\lam)\}\subset
H^0(G/Q,\cl_\lam) $$ given by the path vectors $p_\pi$ as
constructed in \cite{Li$_2$}.

\begin{rem}\rm  The reader not acquainted
with the construction should think of these sections in the following way:
Let $\pi=(\tau_0,\ldots,\tau_r,a_1,\ldots,a_r)$ be an L-S path of shape $\lam$.
Fix $\ell\in\bn$ minimal such that $\ell a_i\in \bn$  for all $i$,
and fix weight vectors $p_{\tau_i}\in H^0(G/Q,\cl_\lam)$ of weight
$-\tau_i(\lam)$ (these weight spaces are one-dimensional). Then $p_\pi$
can be thought of as an algebraic approximation
$$
p_\pi\sim {}^{\ell}\sqrt{p_{\tau_0}^{\ell x_1}p_{\tau_1}^{\ell x_2}p_{\tau_3}^{\ell x_3}\ldots
p_{\tau_r}^{\ell x_{r+1}}},\quad x_i=a_i-a_{i-1}, 1\le i\le r+1.
$$
Note that in the framework of quantum groups at roots of unity the expression above ``makes sense''.
Indeed,  ${}^{\ell} \sqrt{p_{\tau_0}^{\ell x_1}\ldots}$ is then the Frobenius
splitting map at an $\ell$-th root of unity applied to the product, for details see \cite{Li$_2$}.
\end{rem}

The path vectors are $T$-weight vectors of weight $-\pi(1)$. Further, the partial
orders on L-S paths introduced in the section before are closely related to $B$-stable
(respectively $B^-$-stable) submodules spanned by path vectors.

More precisely, we call a subset $S^+\subset B(\lam)$ of L-S paths of shape $\lam$
{\it positive saturated} if for all $\pi,\pi'\in S^+$ the following holds:

$\bullet$ if $\eta\in B(\lam)$ is such that $\pi'>\eta>\pi$, then $\eta\in S^+$.

We say that $S^+\subset B(\lam)$ is {\it maximally positive
saturated} if for $\pi\in S^+$ and $\eta\in B(\lam)$ the relation
$\eta>\pi$ implies $\eta\in S^+$. \vskip 5pt Similarly,  we call a
subset $S^-\subset B(\lam)$ of L-S paths of shape $\lam$ {\it
negative saturated} if for all $\pi,\pi'\in S^-$ the following
holds:

$\bullet$ if $\eta\in B(\lam)$ is such that $\pi'>^r\eta>^r\pi$, then $\eta\in S^-$.

We say that $S^-\subset B(\lam)$ is {\it maximally negative saturated} if
for $\pi\in S^-$ and $\eta\in B(\lam)$ the relation $\pi>^r\eta$ implies $\eta\in S^-$.

Let $S^+$ be a maximally positive saturated subset and let $S^-$
be a maximally negative saturated subset of $B(\lam)$. The corresponding path vectors
span $T$-submodules of $H^0(G/Q,\cl_\lam)$:
$$
M(S^+)=\langle p_\pi\vert\pi\in S^+\rangle
\quad{\rm\  and\ }\quad M(S^-)=\langle
p_\pi\vert\pi\in S^-\rangle
$$

\begin{thm}[\cite{LLM},\cite{Li$_2$}]\label{Bmodule}
The $T$--submodule $M(S^+)$ is $B$--stable
and the $T$--submodule $M(S^-)$ is a $B^-$--stable submodule of $H^0(G/Q,\cl_\lam)$.
\end{thm}

\begin{cor}
Let $S^+\subset B(\lam)$ be a positive saturated subset and let $M(S^+)$
be the $T$-submodule $\langle p_\pi\vert\pi\in S^+\rangle$ of $H^0(G/Q,\cl_\lam)$ spanned by
the corresponding path vectors. Then $M(S^+)$ admits a $B$--module structure which is isomorphic
to a subquotient of a $B$--stable filtration of $H^0(G/Q,\cl_\lam)$.
\end{cor}

{\it Proof\/.} Let $S_1^+=\{\eta\in B(\lam)\mid\exists\,\pi\in S^+:\eta\ge\pi\}$ be the ``closure''
of $S^+$ with respect to $>$, i.e., $S^+_1$ consists of all paths which are greater or equal to
some element of $S^+$. The set $S_1^+$ is a maximally positive saturated subset and hence $M(S_1^+)$ is a
$B$-submodule of $H^0(G/Q,\cl_\lam)$. The set $S_2^+:=S_1^+ - S^+$ is again a maximally positive saturated subset
because: suppose $\pi\in S_2^+$ and $\eta\in B(\lam)$ is such that $\eta>\pi$. By the definition
of $S_1^+$ and $S_2^+$, there exists an element $\pi_1\in S^+$ such that $\pi>\pi_1$, so
$\eta>\pi>\pi_1$ implies $\eta\in S^+_1$. So either $\eta\in S^+_2$ or $\eta\in S^+$. The
latter is not possible because $\eta,\pi_1\in S^+$ would imply $\pi\in S^+$, in contradiction
to the assumption $\pi\in S^+_2$, which finishes the proof of the claim.

Now $M(S_2^+)$ is hence a $B$-submodule, and the quotient $M(S_1^+)/M(S_2^+)$ is a $B$--module
which, as $T$-module, is isomorphic to $M(S^+)$.\endpf

The corresponding version for negative saturated subsets holds also, the proof is left to the reader.
\begin{cor}
Let $S^-\subset B(\lam)$ be a negative saturated subset and let $M(S^-)$
be the $T$-submodule $\langle p_\pi\vert\pi\in S^+\rangle$ of $H^0(G/Q,\cl_\lam)$ spanned by
the corresponding path vectors. Then $M(S^-)$ admits a $B^-$--module structure which is isomorphic
to a subquotient of a $B^-$--stable filtration of $H^0(G/Q,\cl_\lam)$.
\end{cor}

For $\tau\in W/W_\lam$ let $v_\tau\in V(\lam)$ be a weight vector of weight $\tau(\lam)$, $v_\tau$ is
a so-called {\it extremal weight vector}. The $B$-submodule spanned by the orbit $B.v_\tau$
is called the {\it Demazure module} associated to $\tau$ and is denoted $V_\tau(\lam)$. Similarly,
the $B^-$ submodule spanned by the orbit $B^-.v_\tau$
is called the {\it opposite Demazure module} associated to $\tau$ and is denoted $V^\tau(\lam)$.

\begin{thm}[\cite{LLM},\cite{Li$_2$}]\label{demazurebasis}
The path vector basis is compatible with the Demazure submodules, i.e., the restrictions
$\{p_\pi\vert_{V_\tau(\lam)}\mid \pi\in B_\tau(\lam)\}$ form a basis of $V_\tau(\lam)^*$ and the
restrictions of the other path vectors vanish on the submodule.
Similarly, the restrictions $\{p_\pi\vert_{V^\tau(\lam)}\mid \pi\in B^\tau(\lam)\}$ form a basis
of $V^\tau(\lam)^*$ and the restrictions of the other path vectors vanish.
\end{thm}

Let $\bb(\lam)^*=\{u_\pi\in V(\lam)\vert \pi\in B(\lam)\}$ be the basis of $V(\lam)$ {\it dual} to the
path vector basis of $H^0(G/Q,\cl_\lam)$.

\begin{cor}
The vectors $\{u_\pi\vert \pi\in B_\tau(\lam)\}$ form a basis of
the Demazure module $V_\tau(\lam)$, the vectors $\{u_\pi\vert
\pi\in B^\tau(\lam)\}$ form a basis of the opposite Demazure
module $V^\tau(\lam)$, and the vectors $\{u_\pi\vert \pi\in
B_\tau^\sigma(\lam)\}$ form a basis of the intersection
$V_\tau^\sigma(\lam)=V_\tau(\lam)\cap V^\sigma(\lam)$.
\end{cor}

{\it Proof of the corollary\/}: Theorem~\ref{demazurebasis}
implies that $V_\tau(\lam)$ is the subspace of $V(\lam)$
orthogonal to $\langle p_\pi\vert i(\pi)\not\le\tau\rangle$. Hence
$\{u_\pi\vert i(\pi)\le \tau\}\subset V_\tau(\lam)$, and, again by
Theorem~\ref{demazurebasis}, $\langle u_\pi\vert i(\pi)\le
\tau\rangle= V_\tau(\lam)$. The proof for the opposite Demazure
module is similar and is left to the reader. The statement for the
intersection is then an immediate consequence of the fact that the
basis is compatible with Demazure and opposite Demazure modules.
\endpf

\section{Some graded rings and modules}\label{gradedringsandmodules}
Suppose $\lam\in \Lambda_Q^{++}$, so the line bundle $\cl_\lam$ is
very ample on $G/Q$ and we have a corresponding embedding
$G/Q\hookrightarrow \bp(V(\lam))$. Consider the two rings $$
R=\bigoplus_{m\ge 0} H^0(G/Q,\cl_{m\lam})\qquad\hbox{\rm and\
}\qquad k[G/Q]=\bigoplus_{m\ge 0} k[G/Q]_m, $$ where $k[G/Q]$
denotes the homogeneous coordinate ring of the embedding
$G/Q\hookrightarrow \bp(V(\lam))$ with the usual grading. Since
$\lam$ is ample and $G/Q$ is smooth, one knows that
$H^0(G/Q,\cl_{m\lam})=k[G/Q]_m$ for $m\gg 0$. (Actually, by
standard monomial theory \cite{Li$_2$} or Frobenius splitting
\cite{RR}, one knows that they coincide for all $m$, but we do not
need this later.) Correspondingly we denote for a union of
Richardson varieties $X_\uk^\us$ the associated rings by
$R_\uk^\us$ and $k[X_\uk^\us]$, i.e., $$ R_\uk^\us=\bigoplus_{m\ge
0} H^0(X_\uk^\us,\cl_{m\lam}) \quad\hbox{\rm and\ }\quad
k[X_\uk^\us]=\bigoplus_{m\ge 0} k[X_\uk^\us]_m. $$ Let
$\ci_\uk^\us\subset\co_{G/Q}$ be the ideal sheaf and let
$I_\uk^\us\subset k[G/Q]$ be the homogeneous ideal of
$X_\uk^\us\subset G/Q\subset \bp(V(\lam))$. We consider also the
graded $R$-module, respectively the graded $k[G/Q]$-module: $$
J_\uk^\us=\bigoplus_{m\ge 0}
H^0(G/Q,\ci_\uk^\us\otimes\cl_{m\lam}), \quad\hbox{\rm
respectively\ }\quad I_\uk^\us=\bigoplus_{m\ge 0} (I_\uk^\us)_m $$
For Schubert varieties and opposite Schubert varieties we still
use the notation $X_\tau, I_\tau, X^\tau, I^\tau$ etc. instead of
$X_\tau^{id}, I_\tau^{id}, X^\tau_{\overline{w}_0},
I^\tau_{\overline{w}_0}$ etc.

Our aim is to show in the next {\it four sections that the rings $R_\uk^\us$ and
$k[X_\uk^\us]$ and the modules $J_\uk^\us$ and $I_\uk^\us$ coincide and have
a basis by standard monomials}.

\section{Standard monomials and quadratic relations}\label{smtsection}
Let $\lam\in\Lambda_Q^{++}$, we use the same notation as in section~\ref{gradedringsandmodules}.
We analyze the structure of the homogeneous coordinate ring $k[G/Q]$.
We view $k[G/Q]$ as the subalgebra of $R$ generated by
$H^0(G/Q,\cl_{\lam}) =k[G/Q]_1$. A monomial of path vectors
$p_\upi=p_{\pi_1}\cdots p_{\pi_m}$ is called {\it standard} if
the sequence $\upi=(\pi_1,\ldots,\pi_m)$ is standard. The monomial
$p_\upi$ is called {\it standard on a Richardson variety $X_\tau^\sigma$}
if in addition $\tau\ge i(\upi)\ge e(\upi)\ge\sigma$ for the initial and final directions
of $\upi$. The monomial is called {\it standard on a union of Richardson varieties
$X_{\uk}^{\us}$} if it is standard on at least one of the $X_\tau^\sigma\subset X_{\uk}^{\us}$.

The following relations provide an algorithm to express a
non-standard monomial as a linear combination of standard monomials. To be more
precise, we need to introduce the wedge product of two convex subsets of $W/W_\lam$.

Let $\succeq$ be a the fixed total order refining the Bruhat order.
Given two convex subsets $\pi=(\ut,\ua)$ and $\pi'=(\us,\ub)$ of shape $\lam$, let $\uk$ be the
sequence obtained from $\ut$ and $\us$
$$
\{\kappa_0,\kappa_1,\ldots\}=\{\tau_0,\ldots,\tau_r\}\cup\{\sigma_0,\ldots,\sigma_s\}
$$
by writing the elements $\kappa_j$ in strictly decreasing order with respect to $\succ$.

Next let $\uc$ be the strictly increasing sequence of rational
numbers obtained from $\ua$ and $\ub$ as follows: $c_i$ is half
the sum of the weights of the $\tau_j$ and $\sigma_j$ which are
smaller or equal to $\kappa_i$. More precisely: set
$c_0,a_0,b_0:=0$. If $\kappa_{i-1}=\tau_{j-1}$ and not equal to
one of the $\sigma_m$, then $c_i:=c_{i-1}+(a_j-a_{j-1})/2$; if
$\kappa_{i-1}=\sigma_{j-1}$, and not equal to one of the $\tau_m$,
then $c_i:=c_{i-1}+(b_j-b_{j-1})/2$; if
$\kappa_{i-1}=\tau_{j-1}=\sigma_{m-1}$, then
$c_i:=c_{i-1}+((a_j-a_{j-1})+(b_m-b_{m-1}))/2$.

\begin{defn}\rm
The {\it wedge product} $\pi\wedge\pi'$ of two convex subsets of shape $\lam$  is the convex
subset $(\uk,\uc)$ of shape $2\lam$.
\end{defn}

Note that if $\pi,\pi'$ are L-S paths, then $\pi\wedge\pi'$ is in general not an L-S path.
Let now $\pi,\pi'$ be L-S paths of shape $\lam$. We say that $\pi$ and $\pi'$ have the
{\it same support} if the sequence $\uk$ in $\pi\wedge\pi'$ is strictly decreasing with
respect to the Bruhat ordering $>$. One checks easily the following properties:

\begin{lem}\label{wedgeproduct}
\begin{joliste}
\item[{\rm (i)}] The map $(\pi,\pi')\rightarrow \pi\wedge\pi'$
induces a bijection between the set of standard sequences of length two of L-S
paths of shape $\lam$ and the set of L-S paths of shape $2\lam$.
\item[{\rm (ii)}] If $\pi,\pi'$ have the same support, then $\pi\wedge\pi'$
is an L-S path of shape $2\lam$.
\end{joliste}
\end{lem}

\begin{rem}\label{mwedge}\rm
The wedge product of convex sequences generalizes in the obvious
way to $m$-tuples: We define
$(\uk,\uc)=\pi_1\wedge\ldots\wedge\pi_m$ as follows: The sequence
$\uk$ is the union of the Weyl group cosets occurring in the
$\pi_j$, rewritten in strictly decreasing order. Every coset
$\kappa_i$ in $\uk$ occurs in at least one path; set $w_i^k=$ the
weight of the coset $\kappa_i$, if $\kappa_i$ occurs in $\pi_k$,
and set $w_i^k=0$ if $\kappa_i$ does not occur in $\pi_k$. The
rational number $c_i$ in $\uc$ corresponding to the coset
$\kappa_i$ is $$ c_i={1\over m}\sum_{1\le \ell\le i}\ (\sum_{1\le
k\le m} w_\ell^k), $$ the sum over all paths of all weights of the
cosets smaller or equal to $\kappa_i$, divided by $m$. One checks
easily as above: If the $\pi_i$ are L-S paths of shape $\lam$ and
have the same support, then
$\wedge\upi=\pi_1\wedge\ldots\wedge\pi_m$ is an L-S path of shape
$m\lam$. Further, the map $$
\upi=(\pi_1,\ldots,\pi_m)\longrightarrow
\wedge\upi=\pi_1\wedge\ldots\wedge\pi_m $$ induces a bijection
between the standard sequences of length $m$ of L-S paths of shape
$\lam$, and the set of L-S paths of shape $m\lam$. By
construction, one has $i(\upi)=i(\wedge\upi)$ and
$e(\upi)=e(\wedge\upi)$, and also the weights coincide:
$\upi(1)=\wedge\upi(1)$.
\end{rem}

If $\pi\wedge\pi'$ is an L-S path of shape $2\lam$, then {\rm (i)} implies there
exists a standard sequence $(\eta,\eta')$ such that $\eta\wedge\eta'=\pi\wedge\pi'$.
In the following we write just $\pi\wedge\pi'=(\eta,\eta')$ to indicate the corresponding
standard sequence.

We use the notation $(\eta_1,\eta_2)\succeq \pi_1\wedge\pi_2$ if either
$\eta_1\wedge\eta_2\succ \pi_1\wedge\pi_2$, or $\eta_1\wedge\eta_2=\pi_1\wedge\pi_2$
and $\eta_1\succeq \pi_1$. The notation $ \pi_1\wedge\pi_2\succeq^r (\eta_1,\eta_2)$
is defined in the same way: either $\pi_1\wedge\pi_2\succeq^r \eta_1\wedge\eta_2$,
or $\pi_1\wedge\pi_2 = \eta_1\wedge\eta_2$ and $\pi_2\succeq^r\eta_2$.

\begin{thm}[\cite{LLM}]\label{quadraticrelation}
If neither $p_{\pi_1}p_{\pi_2}$ nor $p_{\pi_2}p_{\pi_1}$ is
standard, then there exist standard monomials $p_{\eta_1}p_{\eta_2}\in H^0(G/Q,\cl_{2\lam})$
such that
$$
p_{\pi_1}p_{\pi_2} =\sum a_{\eta_1,\eta_2} p_{\eta_1}p_{\eta_2}
$$
where the coefficient $a_{\eta_1,\eta_2}\not=0$ only if $(\eta_1,\eta_2)\succeq
\pi_1\wedge\pi_2\succeq^r (\eta_1,\eta_2)$. Further, $\pi_1\wedge\pi_2 = \eta_1\wedge\eta_2$
is possible only if $\pi$ and $\pi'$ have the same support, and then $a_{\eta_1,\eta_2}=1$
for the standard sequence $(\eta_1,\eta_2)=\pi_1\wedge\pi_2$.
\end{thm}

Sometimes it is more convenient to formulate the relation as follows:

\begin{cor}\label{quadraticrelationII}
If neither the monomial $p_{\pi_1}p_{\pi_2}$ nor the monomial $p_{\pi_2}p_{\pi_1}$ is
standard, then, in the quadratic relation above, the coefficient of $p_{\eta_1}p_{\eta_2}$
is non-zero only if $\eta_1>\pi_1,\pi_2$ and $\pi_1,\pi_2 >^r\eta_2$.
\end{cor}

{\it Proof.}
This is a consequence of the proof of the theorem above in \cite{LLM}.
The main point we will use from the proof is that $a_{\eta_1,\eta_2}\not=0$ only if
$\eta_1\wedge\eta_2>\pi_1\wedge\pi_2$ in the partial order on convex subsets.
Now on the one hand, the pairs $(\eta_1,\eta_2)$ are standard sequences, so the
wedge product $\eta_1\wedge\eta_2$ is actually independent of the chosen ordering.
On the other hand, $\pi_1\wedge\pi_2$ depends on the choice of the total order on
$W/W_\lam$, so we are free to choose an appropriate total ordering.

Let $\pi_1=(\ut,\ua)$, $\pi_2=(\usigma,\ub)$ and $\eta_1=(\ukappa,\uc)$.
Suppose first $\tau_0$ and $\sigma_0$ are not comparable in the Bruhat order.
We can choose a total order $\succ_1$ such that $\tau_0\succ_1\sigma_0$
and a total order $\succ_2$ such that $\sigma_0\succ_2\tau_0$.
Now $\eta_1\wedge\eta_2>\pi_1\wedge\pi_2$ implies for the first choice
of the total order that $\kappa_0\ge\tau_0$, for the second choice we get as a consequence
$\kappa_0\ge\sigma_0$. Since $\tau_0$, $\sigma_0$ are not comparable,
we must have strict inequality in both cases, in particular $\eta_1> \pi_1,\pi_2$.

If $\tau_0=\sigma_0$, then $\eta_1\wedge\eta_2>\pi_1\wedge\pi_2$ implies by the definition
of the wedge product and the partial order that $\eta_1\ge \pi_1,\pi_2$. So without loss of
generality we can assume in the following $\tau_0>\sigma_0$. Since $\eta_1\wedge\eta_2>
\pi_1\wedge\pi_2$ implies $\kappa_0\ge\tau_0>\sigma_0$, we see that $\eta_1>\pi_2$.
It remains to show: $\eta_1>\pi_1$.

If there exists a $j\le r$ such that $\sigma_0\ge \tau_j$ but
$\sigma_0\not\ge \tau_{j-1}$, then we can choose a total order
such that the sequence of Weyl group cosets in $\pi_1\wedge\pi_2$
is of the form $(\ldots,\tau_{j-1},\sigma_0,\tau_{j},\ldots)$, so
$\eta_1\wedge\eta_2>\pi_1\wedge\pi_2$ implies by the partial
lexicographic ordering that necessarily $\eta_1>\pi_1$. Otherwise
the sequence of Weyl group cosets in $\pi_1\wedge\pi_2$ is of the
form $(\ldots,\tau_{r},\sigma_0,\ldots)$, so
$\eta_1\wedge\eta_2>\pi_1\wedge\pi_2$ implies $\eta_1\ge\pi_1$.
But note that $\eta_1=\pi_1$ and
$\eta_1\wedge\eta_2>\pi_1\wedge\pi_2$ implies $\eta_2\ge\pi_2$ and
hence $\tau_r\ge \sigma_0$, which is not possible by the
assumption that $p_{\pi_1}p_{\pi_2}$ is not standard. It follows
hence also in this case: $\eta_1>\pi_1$.

By replacing $\ge$ and $\succ$ by $\ge^r$ and $\succ^r$, the same arguments show
that $\eta_2>^r\pi_1,\pi_2$.
\endpf

These quadratic relations and the combinatorial character formula are already sufficient to prove:

\begin{prop}\label{homcoordinate}
The homogeneous coordinate ring $k[G/Q]$ has a vector space basis given by the standard
monomials.
\end{prop}

{\it Proof.}
Consider the polynomial algebra $S=k[x_\pi\mid\pi\in B(\lam)]$. We can write any monomial
in $S$ as an ordered product: $x_{\pi_1}^{n_1}\cdots x_{\pi_t}^{n_t}$, where
$\pi_1\succ^r\ldots\succ^r \pi_t$.

We define a monomial order $\succ^r$ on the monomials in S as follows: We say
$x_{\pi}\succeq^r x_{\pi'}$ if $\pi\succeq^r\pi'$. If $m_1,m_2$ are two ordered
monomials, then we say $m_1\succ^r m_2$ if either the degree of $m_1$ is strictly greater
than the degree of $m_2$, or, if the degree of the two coincides, then we say $m_1\succ^r m_2$
if $m_1$ comes before $m_2$ in the reverse lexicographic ordering with respect to $\succeq^r$.

Let $m=x_{\pi_1}\cdots x_{\pi_s}$ be an ordered monomial. By abuse of notation we call $m$ a standard
monomial in $S$ if the sequence $(\pi_1,\ldots,\pi_s)$ is standard. Denote $E\subset S$ the ideal generated
by the type of quadratic relations as in Theorem~\ref{quadraticrelation}, i.e., $E$ is the ideal
generated by the relations
$$
x_{\pi_1}x_{\pi_2} -\sum a_{\eta_1,\eta_2} x_{\eta_1}x_{\eta_2},
$$
where the set of generators runs over all ordered monomials $x_{\pi_1}x_{\pi_2}$ which
are not standard, and the coefficients $a_{\eta_1,\eta_2}$ are as in Theorem~\ref{quadraticrelation}.

For $f\in S$, let $in(f)$ be the initial term, i.e., the greatest
monomial of $f$ with respect to the chosen monomial ordering
$\succ^r$.

Let $m=x_{\mu_1}\cdots x_{\mu_s}\in S$ be an ordered monomial
which is not standard, say $x_{\mu_j}x_{\mu_{j+1}}$ is not
standard. By the relations above, we have,
$x_{\mu_j}x_{\mu_{j+1}}-\sum a_{\eta_1,\eta_2}
x_{\eta_1}x_{\eta_2}$ is an element of $E$, and
$\mu_j\wedge\mu_{j+1}\succ^r (\eta_1,\eta_2)$ for all standard
monomials with nonzero coefficient $a_{\eta_1,\eta_2}$. One checks
easily that this implies $x_{\mu_j}x_{\mu_{j+1}}\succ^r
x_{\eta_1}x_{\eta_2}$ for all standard monomials with nonzero
coefficient $a_{\eta_1,\eta_2}$.

Set $m_1=x_{\mu_1}\cdots x_{\mu_{j-1}}$ and $m_2=x_{\mu_{j+2}}\cdots x_{\mu_s}$, then
$m_1x_{\mu_j}x_{\mu_{j+1}}m_2-\sum a_{\eta_1,\eta_2}m_1x_{\eta_1}x_{\eta_2}m_2$ is an element
of $E$ with initial term $m$. So any non-standard monomial $m\in S$ occurs as the initial
term of an element of $E$.

It follows by Macaulay's Theorem (\cite{Eis}, Theorem 15.3), that the images of the standard
monomials form a generating set for the vector space $S/E$. The canonical epimorphism
$S\rightarrow k[G/Q]$ defined by $x_\pi\mapsto p_\pi$ factors through $S/E$ (Theorem~\ref{quadraticrelation}).
Since $k[G/Q]_m=H^0(G/Q,\cl_{m\lam})\simeq V(m\lam)^*$ for $m\gg 0$, it follows by the
character formula (Theorem~\ref{standardpathseq})
that the number of standard monomials is equal to the dimension of $k[G/Q]_m$, so they
form in fact a basis. It follows that the standard monomials in $(S/E)_m$ are linearly
independent for $m\gg 0$.

Note that $\pi=(id)$ is an L-S path of shape $\lam$. For all
$\ell>0$ and any standard monomial $x_{\pi_1}\cdots x_{\pi_t}$, the monomial
$x_{\pi_1}\cdots x_{\pi_t}x_{id}^\ell$ of degree $t+\ell$ is again standard.
So any linear dependence relation between standard monomials of low degree
can be made into a linear dependence relation between standard monomials of high
degree by multiplying them with a power of $x_{id}$. As a consequence
we see that the standard monomials are linearly independent for all $m\ge 0$,
and the map $S/E\rightarrow k[G/Q]$ is in fact an isomorphism.
\endpf

\section{Ideals and coordinate rings}\label{ideal}

We analyze now the homogeneous coordinate ring and the defining ideal of a union of
Richardson varieties. Throughout this section we assume $\lam\in \Lambda_Q^{++}$
and we consider the embedding
$$
X_\uk^\us\subset G/Q \hookrightarrow \bp(V(\lam)).
$$
\begin{thm}\label{idealtheorem}
The ideal $I_\uk^\us\subset k[G/Q]$ has a basis given by the standard monomials which are
not standard on the union of Richardson varieties $X_\uk^\us$, and the homogeneous
coordinate ring $k[X_\uk^\us]$ has as basis the restrictions of the monomials
standard on $X_\uk^\us$.
Further, the scheme theoretic intersection of two such unions is reduced and is again
a union of Richardson varieties. In particular, a Richardson variety $X_\tau^\kappa$
is non-empty if and only if $\tau\ge \kappa$, and, in this case, $\dim X_\tau^\kappa=\ell(\tau)-\ell(\kappa)$
and it is the scheme theoretic intersection of the Schubert variety $X_\tau$ and the opposite
Schubert variety $X^\kappa$.
\end{thm}

{\it Proof.} We consider first Schubert varieties. By definition,
the linear span of the cone $\hat{X}_\tau$ over $X_\tau\subset  \bp(V(\lam))$ is the Demazure submodule
$V_\tau(\lam)$. If a standard monomial $p_{\pi_1}\cdots p_{\pi_m}$ is not standard on
$X_\tau$, then $p_{\pi_1}$ is not standard on $X_\tau$ and hence
vanishes on $V_\tau(\lam)$ (Theorem~\ref{demazurebasis}). It follows that the standard monomials
which are not standard on $X_\tau$ lie in $I_\tau$.

Since $\cl_\lam$ is ample one knows that $V(m\lam)^*\simeq
k[G/Q]_{m}$ for $m\gg 0$, and the latter has a basis given by
standard monomials of degree $m$. The monomials which are not
standard on $X_\tau$ vanish on $X_\tau$ and hence on
$V_\tau(m\lam)$, the linear span of the cone $\hat{X}_\tau
\hookrightarrow V(m\lam)$. Theorem~\ref{demazurebasis} (applied to
the weight $m\lam$) together with Remark~\ref{mwedge} implies that
$\dim V_\tau(m\lam)^*$ is equal to the number of standard
monomials on $X_\tau$ of length $m$, so the restrictions form a
basis. This proves the linear independence for $m\gg 0$. Now one
can use the same arguments as in the proof of
Proposition~\ref{homcoordinate} to show that the standard
monomials of shape $m\lam$ on $X_\tau$  are linearly independent
for all $m\ge 0$. The proof for an opposite Schubert variety is
the same.

Let  $I_\tau^\kappa$ be the
ideal $I_\tau^\kappa=I_\tau+I^\kappa\subset k[G/Q]$. The vanishing
set of $I_\tau^\kappa$ is $X_\tau^\kappa$. The standard monomials not standard on
$X_\tau^\kappa$ form a basis for $I_\tau^\kappa$. Since
$p_\pi\not\in I_\tau^\kappa$ if and only if $\tau\ge i(\pi)\ge e(\pi)\ge\kappa$, this
shows for $\tau\not\ge\kappa$ that $I_\tau^\kappa=k[G/Q]$ and hence $X_\tau^\kappa=\emptyset$.
We assume in the following: $\tau\ge \kappa$.

It remains to show that $k[G/Q]/I_\tau^\kappa$ has no nilpotent elements, this will follow
from the quadratic relations. Let $f=\sum a_\upi p_\upi\in k[G/Q]$ be a linear combination
of standard monomials on $X_\tau^\kappa$, all of the same degree. Let $\ueta$ be a standard
sequence of L-S path of shape $\lam$ such that $a_\ueta\not=0$, and suppose that $\ueta$
is minimal with respect to ``$\succ$'' with this property. By applying
Theorem~\ref{quadraticrelation} repeatedly, one obtains in $k[G/Q]$:
$$
f^m=a_\ueta^m p_{{\wedge^m\ueta}}+\sum a_\usigma p_{\usigma},
$$
where $\wedge^m\ueta$ is the unique minimal (with respect to ``$\succ$'') element with
non-zero coefficient which occurs in the expression of $f^m$ as a linear combination
of standard monomials. By assumption, $\ueta$ is standard on $X_\tau^\kappa$, so
$\wedge^m\ueta$ is standard on $X_\tau^\kappa$ too. It follows that $f^m\not=0$
in $k[G/Q]/I_\tau^\kappa$ for all $m\ge 1$, and hence $k[G/Q]/I_\tau^\kappa$ has no nilpotent elements.
As an immediate consequence we see that $X_\tau^\kappa$ is the scheme theoretic intersection $X_\tau\cap X^\kappa$.

Let $\tau=\sigma_0>\sigma_1>\ldots>\sigma_r=\kappa$, $r=\ell(\tau)-\ell(\kappa)$ be a strictly
decreasing sequence of elements in $W/W_Q$. Then $\pi_i=(\sigma_i)$ is an L-S path of
shape $\lam$, and the subring $k[p_{\pi_0},\ldots,p_{\pi_r}]\subset k[X_\tau^\kappa]$
is isomorphic to a polynomial ring in $r+1$ variables since all monomials are standard.
There are only a finite number of such strictly decreasing sequences, and these are the only
generators such that products of powers are again standard. So the number of standard monomials of
a fixed degree $m$ is approximately the same as the number of monomials of degree $m$ of a polynomial
ring in $r+1$ variables. So the Hilbert polynomial of $k[X_\tau^\sigma]$ is of degree
$r$ and hence $\dim X_\tau^\sigma=r$.

Next consider a pair of sequences $\uk=(\kappa_1,\ldots,\kappa_t)$
and $\usigma=(\sigma_1,\ldots,\sigma_t)$ such that
$\kappa_j\ge\sigma_j$ for all $j$, and let $X_\uk^\us$ be the
corresponding union of Richardson varieties. Denote
$I_\uk^\us\subset k[G/Q]$ its vanishing ideal, then $p_\upi\in
I_\uk^\us$ for all standard monomials $p_\upi$ not standard on
$X_\uk^\us$. If $f\in I_\uk^\us$, then write $f=\sum a_\upi
p_\upi$ as a linear combination of standard monomials. Since the
restriction of a standard monomial to a Richardson variety either
vanishes or is standard and forms part of a basis, the vanishing
of $f$ on $X_\uk^\us$ implies that $f$ vanishes on all
$X_\kappa^\sigma\subset X_\uk^\us$. Hence none of the $p_\upi$ is
standard on $X_\uk^\us$, so the standard monomials $p_\upi$, not
standard on $X_\uk^\us$, form a basis for $I_\uk^\us$.

The combinatorial description of the ideals shows immediately that the scheme theoretic intersection
$X_\uk^\us\cap X_{\uk'}^{\us'}$ of two such unions is reduced and is again a union of Richardson
varieties.
\endpf
\begin{rem}\rm
The dimension formula and the irreducibility for Richardson
varieties also follow more geometrically from the description of
$Be_\tau\cap B^-e_\kappa$ by Deodhar \cite{Deo}, who showed that
the intersection of the latter is isomorphic to a product of $k$'s
and $k^*$'s, and the  proof by Richardson \cite{Ri}, that the
closure of such an intersection is a Richardson variety; the
reducedness of intersections and unions of Richardson varieties
follow from  \cite{Ram1}. Further, using the methods of
\cite{RR,Ram1,Ram2}, one may also deduce that Richardson varieties
are normal and Cohen-Macaulay (see also \cite{brionlak}, Lemma 1).
As remarked in the Introduction, one of the main aims of this
paper is to present standard monomial theoretic proofs for results
on Richardson varieties and their unions, proof of Theorem
\ref{idealtheorem} being one such example; of course, along the
way, standard monomial theoretic approach also gives us additional
results, for instance, a standard monomial basis for ideals of
Richardson varieties and their unions (as in \ref{idealtheorem}).
\end{rem}

\section{Filtrations of ideals and ideal sheaves}\label{idealsheaves}

The standard monomial bases of the coordinate rings and defining ideals
suggest certain vector space decompositions: For
example, let $X_\tau^\sigma$ be a Richardson variety and suppose $\pi$
is an L-S path of shape $\lam\in \Lambda_Q^{++}$ standard on $X_\tau^\sigma$. The standard monomials of degree
$m$ on $X_\tau^\sigma$ starting with $p_\pi$ are  all of the form $p_\pi p_{\upi'}$,
where $p_{\upi'}$ is a standard monomial of degree $(m-1)$ such that
$$
e(\pi)\ge i(\upi')\ge e(\upi')\ge \sigma.
$$
As a vector space, the space of standard monomials
of degree $m$ on $X_\tau^\sigma$ starting with $p_\pi$ can be identified with the
space of standard monomials of degree $(m-1)$, standard on $X_{e(\pi)}^\sigma$.
The aim of this section is to formulate this vector space decomposition more precisely
in terms of filtrations and associated graded modules and sheaves.

We have two types of boundaries of a Richardson variety $X_\tau^\sigma$, the {\it positive}
and the {\it opposite} or {\it negative} boundary:
$$
\delp X_\tau^\sigma=\bigcup_{\sigma\le\kappa<\tau} X_\kappa^\sigma\qquad\hbox{\rm and\ }\qquad
\delm X_\tau^\sigma=\bigcup_{\sigma<\kappa\le\tau} X_\tau^\kappa
$$
Correspondingly, let $\delp I_\tau^\sigma\subset k[X_\tau^\sigma]$ be the defining ideal
of $\delp X_\tau^\sigma$ and let $\delm I_\tau^\sigma\subset k[X_\tau^\sigma]$ be the defining
ideal of $\delm X_\tau^\sigma$.

The standard monomials (or rather the images), standard on $X_\tau^\sigma$, form a basis
of $k[X_\tau^\sigma]$. By Theorem~\ref{idealtheorem}, $\delp I_\tau^\sigma$ has a basis
given by the standard monomials $p_{\upi}$ on $X_\tau^\sigma$ such
that $i(\upi)=\tau$, and $\delm I_\tau^\sigma$ has a basis given by the standard
monomials $p_{\upi}$ on $X_\tau^\sigma$ such that $e(\upi)=\sigma$.

The partial orders on L-S paths introduced in section~\ref{pathmodel}
are not only strongly related to the $B$--module structure (see section~\ref{pathvector})
but also to the ideal structure of the homogeneous coordinate ring of Richardson varieties.

%
Denote by $S_\tau^+$ {\it the set of L-S paths of shape $\lam$
such that $i(\pi)=\tau$}. Note that this set is positive saturated
(see section~\ref{pathvector}) with a unique maximal element: the
path $(\tau)$.

\begin{thm}\label{idealfiltration}
Let $X_\tau^\sigma$ be a Richardson variety and let $\delp
X_\tau^\sigma$ be its positive boundary. Consider a sequence of
subsets $S_j\subset S_\tau^+$ such that $\vert S_j\vert=j$,
$(\tau)\in S_j$ for $j>0$, and $S_j$ is positive saturated for all
$j$: $$ S_0=\emptyset\subset S_1\subset\ldots \subset
S_{N-1}\subset S_N= S_\tau^+. $$ Set $(\delp
I_\kappa^\sigma)_j=\sum_{\pi\in S_j} k[X_\kappa^\tau]p_{\pi}$,
then the filtration by $T$-stable ideals $$ 0\subset (\delp
I_\kappa^\sigma)_1\subset\ldots\subset (\delp I_\kappa^\sigma)_N
=\delp I_\kappa^\sigma $$ is such that the subquotients are as
$k[X_\kappa^\sigma]$-modules and $T$-modules isomorphic to: $$
(\delp I_\kappa^\sigma)_j/(\delp I_\kappa^\sigma)_{j-1}\simeq
k[X_{e(\pi)}^\sigma](-1)\otimes\chi_{-\pi(1)}, \quad
\{\pi\}=S_j-S_{j-1}, $$ where $k[X_{e(\pi)}^\sigma]$ denotes the
homogeneous coordinate ring of $X_{e(\pi)}^\sigma\subset
X_{\kappa}^\sigma$. The isomorphism is induced by the graded
morphism $$ (\delp I_\kappa^\sigma)_j\rightarrow
k[X_{e(\pi)}^\sigma](-1), \qquad \big(p_\pi f_\pi+\sum_{\eta\in
S_{j-1}} p_{\eta} f_\eta\big)\mapsto
f_\pi\vert_{X_{e(\pi)}^\sigma}. $$
\end{thm}

\begin{rem}\rm
\begin{joliste}
\item[{\rm(i)}] There is a corresponding obvious version of the theorem for the ideal $\delm I_\kappa^\sigma$.
\item[{\rm(ii)}] If $\sigma=id$, i.e., $X_\kappa^\sigma= X_\kappa$ is a Schubert variety, then the filtration
can in addition be chosen to be $B$-equivariant.
\item[{\rm(iii)}] If $\kappa$ is the class of the longest word in $W/W_\lam$, i.e., $X_\kappa^\sigma= X^\sigma$
is an opposite Schubert variety, then the filtration for the ideal $\delm I_\kappa^\sigma$ can in
addition be chosen to be $B^-$-equivariant.
\end{joliste}
\end{rem}

{\it Proof.}
Consider the vector space $J_j\subset k[X_\kappa^\sigma]$ spanned by all standard monomials starting
with a $p_{\pi}$, $\pi\in S_j$, then $J_j\subset (\delp I_\kappa^\sigma)_j$. Since
$S_j$ is positive saturated, Corollary~\ref{quadraticrelationII} implies in fact
that $J_j$ is an ideal and hence, by the definition, $J_j=(\delp I_\kappa^\sigma)_j$. So the latter
has the standard  monomials on $X_\kappa^\sigma$ starting with a $p_{\pi}$, $\pi\in S_j$, as a basis.

This shows that the subquotient $(\delp I_\kappa^\sigma)_j/(\delp I_\kappa^\sigma)_{j-1}$ has as basis
the images of the standard monomials on $X_\kappa^\sigma$ starting with $p_{\pi}$,
$\{\pi\}=S_j-S_{j-1}$.

Let $f=\sum_{\eta\in S_j} p_\eta f_\eta$ be an element of $(\delp I_\kappa^\sigma)_j$. Let
$\pi\in S_j$ be the unique element such that $\pi\not\in S_{j-1}$. Write $f_\pi=\sum a_\ueta p_\ueta
\in k[X_\kappa^\sigma]$ as a linear combination of standard monomials. Let $f_{\pi,1}=\sum a_\ueta p_\ueta$ be
the sum of those summands in the expression above such that $p_{\pi}p_\ueta$ is standard, and set
$f_{\pi,2}=f_\pi-f_{\pi,1}$.

The quadratic relations (Corollary~\ref{quadraticrelationII})
imply that $p_{\pi}f_{\pi,2}\in (\delp I_\kappa^\sigma)_{j-1}$,
and Theorem~\ref{idealtheorem} shows that the restriction of
$f_\pi$ to $X_{e(\pi)}^\sigma$ coincides with the corresponding
restriction of $f_{\pi,1}$. It follows that the map $$ (\delp
I_\kappa^\sigma)_j\rightarrow k[X_{e(\pi)}^\sigma](-1), \qquad
\big(p_\pi f_\pi+\sum_{\eta\in S_{j-1}} p_{\eta}
f_\eta\big)\mapsto f_\pi\vert_{X_{e(\pi)}^\sigma} $$ is a well
defined graded $k[X_\kappa^\sigma]$-module homomorphism. One
checks easily that the map is surjective, has $(\delp
I_\kappa^\sigma)_{j-1}$ as kernel and is $T$-equivariant up to a
twist by the character corresponding to the weight of $p_{\pi}$.
But this is the same as to say that that the homomorphism induces
an isomorphism between the $T$- and $k[X_\kappa^\sigma]$-modules
$(\delp I_\kappa^\sigma)_j/(\delp I_\kappa^\sigma)_{j-1}$ and
$k[X_{e(\pi)}^\sigma](-1)\otimes\chi_{-\pi(1)}$.
\endpf

We consider the corresponding sheaves. Since $\cl_\lam$ is a very ample sheaf on $G/Q$ and the variety is
smooth, one knows that $H^0(G/Q,\cl_{m\lam}) =k[G/Q]_m$ for $m\gg 0$.
Let $\co_{X_\kappa^\sigma}$ be the structure sheaf of $X_\kappa^\sigma\subset G/Q$ and let $\ci_\kappa^\sigma\subset
\co_{G/Q}$ be its sheaf of ideals, i.e., we have the exact sequence
$$
0\rightarrow \ci_\kappa^\sigma\rightarrow\co_{G/Q}\rightarrow \co_{X_\kappa^\sigma}\rightarrow 0.
$$
Again $\cl_\lam$ being ample, for $m\gg 0$ we get an exact sequence:
$$
0\rightarrow H^0(G/Q,\ci_\kappa^\sigma\otimes\cl_{m\lam})\rightarrow H^0({G/Q},\cl_{m\lam})\rightarrow
H^0(X_\kappa^\sigma,\cl_{m\lam})\rightarrow 0
$$
and hence an isomorphism $k[X_\kappa^\sigma]_m\simeq H^0(X_\kappa^\sigma,\cl_{m\lam})$ for $m\gg 0$.
Similarly, let $\delp \ci_\kappa^\sigma$ be the sheaf of ideals of $\delp X_\kappa^\sigma\subset X_\kappa^\sigma$.
The exact sequence
$$
0\rightarrow \delp \ci_\kappa^\sigma \rightarrow\co_{X_\kappa^\sigma}\rightarrow \co_{\delp X_\kappa^\sigma}\rightarrow 0.
$$
induces for $m\gg 0$ an exact sequence:
$$
0\rightarrow H^0(X_\kappa^\sigma, {\delp \ci_\kappa^\sigma}\otimes\cl_{m\lam})\rightarrow
H^0({X_\kappa^\sigma},\cl_{m\lam})\rightarrow  H^0(\delp X_\kappa^\sigma,\cl_{m\lam})\rightarrow 0
$$
and hence isomorphisms $H^0(\delp X_\kappa^\sigma,\cl_{m\lam})\simeq k[\delp X_\kappa^\sigma]_m$ for
$m\gg 0$. For the ideal
$\delp I_\kappa^\sigma$ we get isomorphisms  $H^0(X_\kappa^\sigma, {\delp\ci_\kappa^\sigma}\otimes\cl_{m\lam})\simeq
[\delp I_\kappa^\sigma]_m$ for $m\gg 0$, where the grading is defined by
$[\delp I_\kappa^\sigma]_m=\delp I_\kappa^\sigma \cap k[X_\kappa^\sigma]_m$.

Set $R=\bigoplus_{m\ge 0} H^0(G/Q,\cl_{m\lam})$. Since $\cl_\lam$
is ample, coherent sheaves on $G/Q$ correspond to finitely
generated graded $R$-modules. If $M$ is such a module, the sheaf
is given locally by quotients (in the localization of $M$) of
homogeneous elements of $M$ and $R$ of the same degree.
Recall that graded modules which differ only in a finite number of graded components
give the same coherent sheaf. It follows that the filtration of $\delp I_\kappa^\sigma$
(Theorem~\ref{idealfiltration}) provides a filtration of $\delp\ci_\kappa^\sigma$
as $\co_{G/Q}$-module and as a $T$-sheaf:

\begin{cor}[Pieri-Chevalley type formula]\label{idealsheaffiltration}
The filtration of the vanishing ideal $\delp I_\kappa^\sigma$ induces a filtration of $\delp\ci_\kappa^\sigma$
such that the associated graded sheaf is, as $T$-equivariant sheaf of $\co_{G/Q}$-modules,
the direct sum
$$
\grad \delp\ci_\kappa^\sigma=\bigoplus \co_{X_{e(\pi)}^\sigma}(-1)\otimes \chi_{-\pi(1)}
$$
where the sum runs over all L-S paths $\pi$ of shape $\lam$, standard on $X_\kappa^\sigma$
and such that $i(\pi)=\kappa$.
\end{cor}

\section{Vanishing theorems and standard monomials}\label{smtandvan}
As before, we assume $\lam\in \Lambda_Q^{++}$.

\begin{thm}\label{vanishandbasis}
Let $X_\uk^\us$ be a union of Richardson varieties.
\begin{joliste}
\item[{\rm (i)}] $H^i(X_\uk^\us,\cl_{m\lam})=0$ for all $i\ge 1$ and all $m\ge 1$,
and for an irreducible variety one has in addition $H^i(X_\kappa^\sigma,\co_{X_\kappa^\sigma})=0$
for all $i\ge 1$.
\item[{\rm (ii)}] The standard monomials, standard on $X_\uk^\us$ of degree $m$, form a basis
for $H^0(X_\uk^\us,\cl_{m\lam})$ for all $m\ge 1$.
\end{joliste}
\end{thm}

{\it Proof\/}
The proof is by induction on the maximal dimension of the irreducible components and on
the number of irreducible components of maximal dimension and on $m$.
The theorem holds if $\dim X_\uk^\us=0$.

Assume that the theorem holds for all unions of Richardson varieties
of dimension smaller than $n$.
Suppose now $\dim X_\kappa^\sigma=n$.
Consider the exact sequence
$$
0\rightarrow \delp \ci_\kappa^\sigma\otimes\cl_\lam \rightarrow\co_{X_\kappa^\sigma}\otimes\cl_\lam
\rightarrow \co_{\delp X_\kappa^\sigma}\otimes\cl_\lam\rightarrow 0.
$$
Since $\delp X_\kappa^\sigma$ is a union of Richardson varieties of smaller dimension,
the vanishing theorem holds and hence $H^i(\delp X_\kappa^\sigma,\cl_\lam)=0$ for $i\ge 1$.
Moreover, by induction, the global sections on $\delp X_\kappa^\sigma$ can
be lifted to global sections on $G/Q$, so the restriction map $H^0(X_\kappa^\sigma,\cl_\lam)
\rightarrow H^0(\delp X_\kappa^\sigma,\cl_\lam) $ is surjective.

By Corollary~\ref{idealsheaffiltration}, the sheaf $\delp
\ci_\kappa^\sigma\otimes\cl_\lam$ admits a filtration such that
the associated graded is just a direct sum of structure sheafs of
Richardson varieties of the form $\co_{X_{\kappa'}^\sigma}$, where
$\kappa\ge \kappa'$. By induction, for $\kappa'<\kappa$ the
vanishing for the higher cohomology groups holds. The associated
graded has only one subquotient isomorphic to
$\co_{X_{\kappa}^\sigma}$, the contribution coming from the path
$(\kappa)$. So the long exact sequence in cohomology splits into
isomorphisms $H^i(X_\kappa^\sigma,\co_{X_{\kappa}^\sigma}) \simeq
H^i(X_\kappa^\sigma,\cl_\lam)$ for $i\ge 1$ and a short exact
sequence $$ 0\rightarrow H^0(X_\kappa^\sigma,\delp
\ci_\kappa^\sigma\otimes\cl_\lam) \rightarrow
H^0(X_\kappa^\sigma,\cl_\lam) \rightarrow H^0(\delp
X_\kappa^\sigma,\cl_\lam)\rightarrow 0, $$ so $\dim
H^0(X_\kappa^\sigma,\cl_\lam)=\dim H^0(\delp
X_\kappa^\sigma,\cl_\lam)+ \dim H^0(X_\kappa^\sigma,\delp
\ci_\kappa^\sigma\otimes\cl_\lam)$. Consider the right hand side
of the equation. By induction, the first term is equal to the
number of L-S paths of shape $\lam$ standard on $\delp
X_\kappa^\sigma$, and the second is by the filtration equal to the
number of L-S paths of shape $\lam$ standard on $X_\kappa^\sigma$
and such that $i(\pi)=\kappa$. This together is the number of L-S
paths of shape $\lam$ standard on $X_\kappa^\sigma$. Since
$k[X_\kappa^\sigma]_1\hookrightarrow
H^0(X_\kappa^\sigma,\cl_\lam)$, this proves that the path vectors
form a basis of $H^0(X_\kappa^\sigma,\cl_\lam)$.

For $m>1$ one gets an exact sequence: $$ 0\rightarrow \delp
\ci_\kappa^\sigma\otimes\cl_{m\lam}
\rightarrow\co_{X_\kappa^\sigma}\otimes\cl_{m\lam} \rightarrow
\co_{\delp X_\kappa^\sigma}\otimes\cl_{m\lam}\rightarrow 0. $$ The
same arguments as above (induction on the dimension or on $m$ and
the filtration of $\delp \ci_\kappa^\sigma\otimes\cl_{m\lam}$)
show that the associated long exact sequence in cohomology splits
into isomorphisms $H^i(X_\kappa^\sigma,\cl_{(m-1)\lam}) \simeq
H^i(X_\kappa^\sigma,\cl_{m\lam})$ for $i\ge 1$ and a short exact
sequence for the global sections. The same counting argument as
above proves that the standard monomials form a basis of
$H^0(X_\kappa^\sigma,\cl_{m\lam})$. Further, for $i\ge 1$ we have
isomorphisms: $$
H^i(X_\kappa^\sigma,\co_{X_{\kappa}^\sigma})\simeq
H^i(X_\kappa^\sigma,\cl_{\lam}) \simeq
H^i(X_\kappa^\sigma,\cl_{m\lam}),\ \forall m\ge 1 $$ Since
$\cl_\lam$ is ample, it follows that
$H^i(X_\kappa^\sigma,\cl_{m\lam})=0$ for all $i\ge 1, m\ge 0$.

Suppose now $X_\uk^\us$ is a union of Richardson varieties. Let
$X_\kappa^\sigma$ be an irreducible component of maximal dimension
and denote $X_{\uk'}^{\us'}$ the union of the remaining
irreducible components. We have an exact sequence: $$ 0\rightarrow
\co_{X_\uk^\us}\rightarrow \co_{X_\kappa^\sigma}\oplus
\co_{X_{\uk'}^{\us'}}\rightarrow \co_{X_{\uk'}^{\us'}\cap
X_\kappa^\sigma}\rightarrow 0. $$ If we tensor the sequence with
$\cl_{m\lam}$, $m\ge 1$, then we know by induction on the
dimension, respectively the number of irreducible components of
maximal dimension that the higher cohomology groups of the second
and the last term vanish, and further, for the global sections the
last map is surjective. It follows that the higher cohomology
groups $H^i(X_\uk^\us,\cl_{m\lam})$ vanish.

It remains to count the dimensions of the spaces of global
sections. Since one has standard monomial theory for the last two
terms by induction on the dimension, respectively the number of
irreducible components of maximal dimension, one checks easily
that in fact $\dim H^0(X_\uk^\us,\cl_{m\lam})$ is equal to the
number of standard monomials on $X_\uk^\us$ of degree $m$. The
same argument as above proves hence that the standard monomials
form a basis of $\dim H^0(X_\uk^\us,\cl_{m\lam})$.
\endpf
\begin{rem}\rm
The result (i) in Theorem \ref{vanishandbasis} is also proved in
\cite{brionlak} (cf. \cite{brionlak}, Lemma 5).
\end{rem}

\section{The non-regular case and pointed unions}\label{nonregular}
Fix a dominant weight $\lam\in \Lam_Q^+$. Note that in the general
case not all results carry over. We have the following two
examples (cf. \cite{brionlak}, Remark following Lemma 5).

Consider $G=SL_2$, so $G/B=\bp^1$. Let $X_\uk^\us=\{0,\infty\}$ be
the union of the two $T$-fixed points. Then $\dim
H^0(G/B,\co_{G/B})=1$ but $\dim H^0(X_\uk^\us,\co_{X_\uk^\us})=2$,
so the restriction map for global sections is not surjective.

 In $\bp^1\times\bp^1$, let $X_\uk^\us$ be the union of the four
$\bp^1$'s: $\{0\}\times \bp^1$, $\{\infty\}\times \bp^1$,
$\bp^1\times \{0\}$ and $\bp^1\times \{\infty\}$, then
$H^1(X_\uk^\us,\co_{X_\uk^\us})\not=0$.

But we get a SMT (standard monomial theory) and vanishing theorems also in the non-regular case
for a special class of unions of Richardson varieties. In the inductive procedure used to
construct standard monomial theory, we need a class of varieties which
includes {\it i)} all Richardson varieties, {\it ii)} their boundaries $\delp X_\tau^\kappa$ and
$\delm X_\tau^\kappa$, and {\it iii)} if $Y$ is in this class and $Y=X\cup Y'$ is such that $X$ is irreducible,
then $X$, $Y\cap X$ and $Y'$ are also in this class.

The example above shows: we can expect that the global restriction map is surjective for all varieties in
this class only if it is not possible to construct a non-trivial union of points by using the operations in
{\it ii)} and {\it iii)}.
This leads to the following definition:

\begin{defn}\label{pointedRichardson}
{\it A union of Richardson varieties $Y$ is called pointed\/} if there exists a $\kappa\in W/W_Q$
such that $Y=\bigcup_j X_\kappa^{\sigma_j}$, or if there exists a
$\sigma\in W/W_Q$ such that $Y=\bigcup_j X_{\kappa_j}^\sigma$.
\end{defn}

Simple examples for such pointed unions are arbitrary unions of Schubert varieties or
arbitrary unions of opposite Schubert varieties. It is easy to check that the
class of pointed unions of Richardson varieties has the three properties above.

A first step to construct SMT is the generalization of
Corollary~\ref{idealsheaffiltration}, but first the definition of
the boundary  has to be adjusted. For $\tau\in W/W_Q$, set
$\tbar\equiv\tau \bmod W_\lam$, and denote $p_\tbar$ the extremal
weight vector in $H^0(G/Q,\cl_\lam)$ corresponding to the L-S path
$(\tbar)$ of type $\lam$.

\begin{prop}\label{emptylambdaboundary}
Let $P_\lam\supset Q$ be the parabolic subgroup associated
to the weight $\lam$, and for $\tau\ge \kappa\in W/W_Q$ consider the projection:
$$
\begin{array}{rccc}
\phi:& X_\tau^\kappa&\longrightarrow &G/P_\lam\\
  &\hskip 15pt  \searrow   &                & \nearrow\hskip 15pt\\
  &              &    G/Q   & \\
\end{array}
$$ Either $X_\tau^\kappa\subset \phi^{-1}(e_\tbar)$, in which case
the line bundle $\cl_\lam\vert_{X_\tau^\kappa}$ is trivial, or the
preimage of set of zeros of $p_\tbar$ below, called the boundary
of $X_\tau^\kappa$ with respect to $\lam$ (or just the
$\lambda$-{\it boundary}), is not empty: $$ \delp_\lam
X_\tau^\kappa = \{y\in X_\tau^\kappa\mid
p_\tbar(y)=0\}\not=\emptyset $$
\end{prop}

\begin{rem}\rm
 We use the notation $ \delp_\lam X_\tau^\kappa$ in the following for the corresponding
 variety with its induced reduced structure.
 \end{rem}

{\it Proof.\/}
Let $Z=\phi(X_\tau^\kappa)$ be the image of the $T$-equivariant map $\phi$.
If $Z=e_\tbar$, then $X_\tau^\kappa\subset \phi^{-1}(e_\tbar)$ and the
restriction  of $\cl_\lam$ to the fibre is constant. Otherwise $\dim Z\ge 1$ and hence $Z$ admits
at least two $T$-fixed points, say  $e_\tbar$ and $e_\sbar$ for some $\sbar\not=\tbar$.
Let $y\in X_\tau^\kappa$ be such that $\phi(y)=e_\sbar$, then $p_\tbar(y)=p_\tbar(e_\sbar)=0$
and hence $y\in \delp_\lam X_\tau^\kappa$, so the latter is not empty.
\endpf

\begin{rem}\label{boundarydescription}\rm
Let $\tbar$ be as above and let $\delp X_\tbar=\bigcup
X_{\sbar_i}\subset G/P_\lam$ be the boundary of $X_\tbar$. Let
$X_{\sigma_i}\subset G/Q$ be the preimage of $X_{\sbar_i}$, then
$$ \delp_\lam X_\tau^\kappa = \bigcup (X_\tau^\kappa\cap
X_{\sigma_i}), $$ and this intersection holds scheme
theoretically. In particular, the boundary is either empty or a
pointed union of Richardson varieties. In view of Proposition
\ref{emptylambdaboundary}, we shall henceforth consider only those
$X_\tau^\kappa$ for which the $\lambda$-boundary is non-empty.
\end{rem}

It  remains to define {\it standardness on $X_\tau^\kappa$\/} for the non-regular case:
If $\pi=(\sigma_0,\ldots,\sigma_p,a_1,\ldots,a_p)$ is an L-S path of shape $\lam$, then
$\pi$ is called {\it standard on\/ $X_\tau^\kappa$\/} if there exist {\it lifts}
$\ts_0,\ldots,\ts_p\in W/W_Q$, i.e., $\ts_i\equiv \sigma_i\bmod W_\lam$ for
all $i=0,\ldots,r$,
such that
$$
\tau\ge\ts_0\ge\ts_1\ge\ldots\ge\ts_p\ge\kappa.
$$
Such a sequence $(\ts_0,\ldots,\ts_p)$ is called a {\it defining chain for $\pi$} (on $X_\tau^\kappa$).

For the $\lambda$-boundary $\delp_\lam X_\tau^\kappa$ of the
Richardson variety $X_\tau^\kappa$, let $\ci_{\delp_\lam
X_\tau^\kappa}$ be the corresponding sheaf of ideals.

\begin{thm}\label{idealfiltrationtwo}
The twisted sheaf $\ci_{\delp_\lam X_\tau}(\lam)=\ci_{\delp_\lam
X_\tau^\kappa}\otimes\cl_\lam$ admits a filtration, such that the
associated graded is isomorphic (as $\co_{X_\kappa^\sigma}$--sheaf
and as $T$-sheaf) to a direct sum of structure sheaves, twisted by
a character: $$ \grad \ci_{\delp_\lam X_\tau^\kappa}(\lam)
\simeq\bigoplus_{\pi\in S}
\co_{X_{e(\pi)}^\kappa}\otimes\chi_{-\pi(1)}, $$ where $S$ is the
set of all L--S paths $\pi$ of shape $\lam$, standard on
$X_\tau^\kappa$ and such that $i(\pi)\equiv\tau\bmod W_\lam$.
\end{thm}

The rest of this section is devoted to the proof of Theorem
\ref{idealfiltrationtwo}. For the rest of this section we fix a
regular dominant weight $\rho\in\Lambda^{++}_Q$. Let $\pi\in
B(\lam)$ be as above and let
$\eta_1,\ldots,\eta_s,\eta_1',\ldots,\eta_q'\in B(\rho)$.

\begin{defn}\rm
A sequence $\ueta=(\eta_1,\ldots,\eta_s,\pi,\eta_1',\ldots,\eta_q')$
of shape $(s\rho,\lam,q\rho)$, is called {\it standard} if there
exists a defining chain $(\ts_0,\ldots,\ts_p)$ for $\pi$ such that
$$
e(\eta_1)\ge i(\eta_2)\ge\ldots\ge e(\eta_s)\ge \ts_0\ge\ldots\ge\ts_p\ge i(\eta'_{1})\ge e(\eta'_1)
\ge \ldots\ge i(\eta'_q).
$$
The sequence $\ueta$ is called {\it standard on $X_\tau^\kappa$} if in addition $\tau\ge i(\eta_1)$
(respectively $\tau\ge \ts_0$ if $s=0$) and $e(\eta'_q)\ge\kappa$ (respectively $\ts_p\ge \kappa$ if $q=0$).
The sequence $(\ts_0,\ldots,\ts_p)$ is then called a {\it defining chain} for
$\ueta$ on $X_\tau^\kappa$.
\end{defn}

The monomial $p_\ueta=p_{\eta_1}\cdots p_{\eta_s} p_\pi p_{\eta'_{1}}\cdots p_{\eta_q'}$
is called {\it standard} (on  $X_\tau^\kappa$) if the sequence $\ueta$ is so.
The monomial is called {\it standard on a union of Richardson varieties} if it is standard on at least one
irreducible component. We have the following global basis given by standard monomials:

\begin{thm}[\cite{Li$_2$}]\label{globalmonom}
The standard monomials of shape $(s\rho,\lam,q\rho)$ form a basis for $H^0(G/Q,\cl_{\lam+(s+q)\rho})$
for all $s,q\ge 0$.
\end{thm}

The notion of the final and initial direction has to be adapted to
standard sequences. Let $\ueta$ be a standard sequence on
$X_\tau^\kappa$. The set of all defining chains for $\ueta$ on
$X_\tau^\kappa$ is partially ordered: We say
$(\ts_0,\ldots,\ts_p)\ge (\ts_0',\ldots,\ts_p')$ if and only if
$\ts_0\ge\ts_0'$, $\ldots$, $\ts_p\ge \ts_p'$. One knows by
Deodhar's Lemma that for this partial order there exists a unique
minimal and a unique maximal defining chain (see for example
\cite{$G/P$-IV},  \cite{$G/P$-V} or \cite{Li$_3$}) .

\begin{defn}\label{nonregularlastchain}\rm
If $s\ge 1$, then set $i(\ueta)=i(\eta_1)$ as before. If $s=0$,
then set $i(\ueta)=\ts_0$ for the  unique minimal defining chain
$(\ts_0,\ldots,\ts_p)$.

If $q\ge 1$, then set $e(\ueta)=e(\eta'_q)$. If $q=0$, then set
$e(\ueta)=\ts_p$ for the unique maximal defining chain
$(\ts_0,\ldots,\ts_p)$.
\end{defn}

\begin{rem}\label{definingandclassic}\rm
If $q\ge 1$ or $s\ge 1$, then $i(\ueta)$ and $e(\ueta)$ are
independent of the Richardson variety $X_\tau^\kappa$. If $q=s=0$,
this is not the case: For $G=SL_3$, $Q=B$, $\lam=\omega_1$,
$\eta=(id)$ and $X_\tau^\kappa=X_{id}^{s_1s_2}$, we have
$i(\eta)=(id)$, but for $X_\tau^\kappa=X_{s_1}^{s_1s_2}$, we have
$i(\eta)=(s_1)$.
\end{rem}

\begin{rem}\label{defininganddirect}\rm
The path model theory provides for all $s,q\ge 0$, $q+s>0$, a natural weight preserving bijection
between the set of all L--S paths $\pi$ of shape $\lam+(s+q)\rho$ such that $i(\pi)=\tau$ and
$e(\pi)=\kappa$ and the set of standard sequences $\ueta$ of shape $(s\rho,\lam,q\rho)$ such that
$i(\ueta)=\tau$ and $e(\ueta)=\kappa$, see \cite{Li$_3$}.
\end{rem}

\begin{prop}\label{smtonenonreg}
The standard monomials of shape $(0,\lam,0)$ are linearly
independent on a pointed union of Richardson varieties $Y$.
Further, for all $s+q\ge 1$, the standard monomials on $Y$ of shape
$(s\rho,\lam,q\rho)$ form a basis for $H^0(Y,\cl_{\lam+(s+q)\rho})$
and the corresponding higher cohomology groups vanish.
\end{prop}

{\it Proof.\/} Theorem~\ref{vanishandbasis} implies for $q+s\ge 1$ that the higher
cohomology groups vanish, and,  by Remark~\ref{defininganddirect},
the number of standard monomials, standard on $Y$, is equal
to the dimension of $H^0(Y,\cl_{\lam+(s+q)\rho})$. So to prove the proposition in this
case, it is sufficient to show that the standard monomials span the space of global
sections. We consider only the case $q\ge1$, the case $s\ge 1$ can be proved similarly.

The proof for a union of (opposite) Schubert varieties $X_\tau$
has been given in \cite{Li$_2$}. Suppose $Y=\bigcup_j
X_{\tau_j}^\kappa$ is a pointed union of Richardson varieties.
Denote by $Z$ the union of Schubert varieties $Z=\bigcup_j
X_{\tau_j}$. The restriction map
\begin{equation}\label{restrictionmapproposition}
H^0(Z,\cl_{\lam+(q+s)\rho})\rightarrow H^0(Y,\cl_{\lam+(q+s)\rho}),
\end{equation}
is surjective by Theorem~\ref{vanishandbasis}. The standard monomials standard on $Y$ are also
standard on $Z$. But a standard monomial $p_\ueta$, standard on $Z$, is not standard on $Y$ if
and only if $e(\ueta)\not\ge\kappa$. But in this case the restriction of $p_\ueta\vert_Y$ vanishes
($q>0$!), so, by the surjectivity of $(\ref{restrictionmapproposition})$, the standard monomials,
standard on $Y$, span the space of global sections.

If $Y=\bigcup_j X_{\tau}^{\kappa_j}$, then let $Z$ be the Schubert variety $Z=X_{\tau}$. As above,
the restriction map on the global sections is surjective, and the same arguments show that a standard
monomial, standard on $Z$ but not on $Y$, vanishes identically on $Y$, which finishes the proof
also in this case.

It remains to consider the case $q=s=0$. Suppose first $Y=\bigcup_j X_{\tau_j}^\kappa$ is a pointed
union of Richardson varieties and $p_{\pi_1},\ldots,p_{\pi_t}\in H^0(G/Q,\cl_\lam)$ are standard on $Y$,
but $(\sum_{i=1}^t a_i  p_{\pi_i})\vert_Y=0$. The sequence $(\kappa)$ is an L-S path of shape $\rho$, let
$p_\kappa$ be the corresponding path vector. By the definition of standardness, the monomials
$p_{\pi_i}p_\kappa$ of shape $(0,\lam,\rho)$ are standard on $Y$.

It follows that the sum $\sum_{i=1}^t a_i p_{\pi_i}p_\kappa=0$
would be a linear dependence relation of standard monomials of shape $(0,\lam,\rho)$, in contradiction
to what has been proved above. For $Y=\bigcup_j X_{\tau}^{\kappa_j}$ the arguments are similar, instead
of $p_\kappa$ one takes the path vector $p_\tau\in H^0(G/Q,\cl_\rho)$ and deduces the linear
independence of the $p_{\pi_1},\ldots,p_{\pi_t}$ by the linear independence of
the $p_\tau p_{\pi_1},\ldots,p_\tau p_{\pi_t}$.
\endpf

We come now to the proof of Theorem~\ref{idealfiltrationtwo}.
\par\noindent
{\it Proof.\/}
The exact sequence $0\rightarrow \ci_{\delp_\lam X_\tau^\kappa}\rightarrow
\co_{X_\tau^\kappa}\rightarrow \co_{\delp_\lam X_\tau^\kappa}\rightarrow 0$ induces an
exact sequence
$$
0\rightarrow \ci_{\delp_\lam X_\tau^\kappa}(\lam)\rightarrow \co_{X_\tau^\kappa}\otimes\cl_{\lam}
\rightarrow \co_{\delp_\lam X_\tau^\kappa}\otimes\cl_{\lam}\rightarrow 0.
$$
Since $\rho$ is regular, one gets for $m\gg 0$ an induced exact sequence
\begin{equation}\label{shortexactnonregular}
\begin{array}{rll}
0\rightarrow H^0(X_\tau^\kappa, \ci_{\delp_\lam X_\tau^\kappa}(\lam)\otimes \cl_{m\rho})\rightarrow &&\\
H^0(X_\tau^\kappa, \cl_{\lam+m\rho})&\rightarrow&
H^0(\delp_\lam X_\tau^\kappa,\cl_{\lam+m\rho})\rightarrow 0.
\end{array}
\end{equation}
Denote $S_\tbar^+$ {\it the set of L-S paths of shape $\lam$,
standard on $X_\tau^\kappa$ and such that $\pi=(\tbar,\ldots)$}.
This set is positive saturated, i.e., if $\pi,\pi'$ are elements
of this set and $\pi'\ge\eta\ge\pi$, then $\eta$ is an element of
this set too (see section~\ref{pathvector}). This set has a unique
maximal element: the path $(\tbar)$. By deleting step by step one
element, it is possible to get a sequence of subsets $$
S_0=\emptyset\subset S_1\subset\ldots \subset S_{N-1}\subset S_N=
S_\tbar^+. $$ such that $(\tbar)\in S_j$ for all $j>0$ and the
$S_j$ are positive saturated for all $j$. Fix $M\gg 0$ so that
$(\ref{shortexactnonregular})$ holds for all $m\ge M$ and consider
the $R$--module $$ I_\tbar^\kappa(\lam):=\bigoplus_{m\ge M}
H^0(X_\tau^\kappa,\cl_{m\rho} \otimes \ci_{\delp_\lam
X_\tau^\kappa}(\lam)). $$ The path vectors $p_\pi$, $\pi\in
S^+_\tbar$, are sections in $H^0(X_\tau^\kappa,\cl_\lam)$ which
vanish on $\delp_\lam X_\tau^\kappa$, this follows immediately
from the description of the boundary in
Remark~\ref{boundarydescription} and the fact that the $p_\pi$
vanishes on $\delp X_\tbar$. Set $$
(I_\tbar^\kappa(\lam))_j=\sum_{\pi\in S_j}\left(\sum_{m\ge
M}H^0(X_\kappa^\tau,\cl_{m\rho}) p_{\pi}\right), $$ then
$(I_\tbar^\kappa(\lam))_j$ is a $T$-stable $R$--submodule of
$I_\tbar^\kappa(\lam)$ and hence this defines a filtration of
$I_\tbar^\kappa(\lam)$ by $T$-stable $R$--submodules $$ 0\subset
(I_\tbar^\kappa(\lam))_1\subset\ldots\subset
(I_\tbar^\kappa(\lam))_N =I_\tbar^\kappa(\lam). $$ The middle term
in the short exact sequence~{(\ref{shortexactnonregular})} has as
basis the standard monomials of shape $(0,\lam,m\rho)$, standard
on $X_\tau^\kappa$, and the last term has as basis the standard
monomials of shape $(0,\lam,m\rho)$, standard on $\delp_\lam
X_\tau^\kappa$ (see Proposition~\ref{smtonenonreg}). Let
$I(\lam+m\rho)$ be the set of standard monomials of shape
$(0,\lam,m\rho)$, standard on $X_\tau^\kappa$ but not standard on
$\delp_\lam X_\tau^\kappa$. These are the standard monomials of
the form $p_\ueta=p_\pi p_{\ueta'}$, $\pi\in S_\tbar^+$ and
$p_{\ueta'}$ of shape $m\rho$. These monomials vanish on
$\delp_\lam X_\tau^\kappa$. It follows that the $p_\ueta$,
$\ueta\in I(\lam+m\rho)$, form a basis of the first term.

The next step is to prove that the basis is compatible with the
filtration. So suppose $p_\pi\in H^0(X_\tau^\kappa,\cl_\lam)$ is
such that $\pi\in S_\tbar^+$ and let $p_\eta$ be a path vector in
$H^0(X_\tau^\kappa,\cl_\rho)$ and standard on $X_\tau^\kappa$.
Suppose the monomial $p_\pi p_\eta$ of shape $(0,\lam,\rho)$ is
not standard on $X_\tau^\kappa$. Since $\rho\in\Lam^{++}$, this
implies automatically that the product is not standard on the
Schubert variety $X_\tau$. Now the quadratic relations in
\cite{LiSe}, Theorem~3.13, show: $$ p_\pi p_\eta =
\sum_{\pi',\eta'} b_{\pi',\eta'} p_{\pi'}p_{\eta'} $$ gives a
presentation of the product as linear combination of standard
monomials, standard on $X_\tau$, where
$(\pi',\eta')\ge(\pi,\eta)\ge^r (\pi',\eta')$. Such a standard
monomial $p_\pi p_\eta$ is then not standard on $X_\tau^\kappa$ if
and only if $e(\eta')\not\ge\kappa$, in which case the monomial
vanishes on $X_\tau^\kappa$. But this implies that we have the
same type of relation for non-standard products also on
$X_\tau^\kappa$ and hence:
\begin{equation}
p_\pi p_\eta = \sum_{(\pi',\eta')\ge(\pi,\eta)\ge^r (\pi',\eta')} b_{\pi',\eta} p_{\pi'}p_{\eta'},
\hbox{\rm\ all\ $p_{\pi'}p_{\eta'}$\ are standard on\ }X_\tau^\kappa
\end{equation}
The sets $S_j$ are positively saturated, and hence (as in the regular case)
the subspace spanned by elements of degree $m$ in the $j$--th filtration part:
$$
(I_{\tbar,m}^\kappa(\lam))_j=H^0(X_{\tau}^\kappa, \ci_{\delp_\lam X_\tau^\kappa}(\lam)\otimes \cl_{m\rho})_j
$$
has a basis given by the standard monomials, standard on $X_\tau^\kappa$, starting with a $p_\pi$,
$\pi\in S_j$. Let $\{\pi\}= S_j-S_{j-1}$. As in the regular case, the map
$$
(I_{\tbar,m}^\kappa(\lam))_j\rightarrow H^0(X_{e(\pi)}^\kappa,\cl_{m\rho}) \qquad
\left(p_\pi f_\pi+\sum_{\eta\in S_{j-1}} p_{\eta} f_\eta\right)\mapsto f_\pi\vert_{X_{e(\pi)}^\kappa}.
$$
induces an isomorphism $(I_{\tbar,m}^\kappa(\lam))_j/(I_{\tbar,m}^\kappa(\lam))_{j-1}\longrightarrow
H^0(X_{e(\pi)}^\tau,\cl_{m\rho})$. Since this holds for all $m\ge M$, this induces
a filtration of the sheaf
$$
0\subset (\ci_{\delp_\lam X_\tau^\kappa}(\lam))_1\subset (\ci_{\delp_\lam X_\tau^\kappa}(\lam))_2\subset
\ldots\subset (\ci_{\delp_\lam X_\tau^\kappa}(\lam))_N=\ci_{\delp_\lam X_\tau^\kappa}(\lam)
$$
such that the associated graded is isomorphic to the direct sum
$$
\grad \ci_{\delp_\lam X_\tau^\kappa}(\lam)=\bigoplus_{\pi\in S} \co_{X_{e(\pi)}^\kappa}\otimes\chi_{-\pi(1)}.
$$
\endpf

\section{SMT in the non-regular case}\label{SMTnonreg}
We have the following generalization of Theorem~\ref{vanishandbasis} to the non-regular case:

\begin{thm}\label{vanishandbasisnonreg}
Suppose $\lam\in\Lambda^+_Q$ is a dominant weight and let $Y$ be a pointed union of Richardson varieties.
\begin{joliste}
\item[{\rm (i)}] $H^i(Y,\cl_{m\lam})=0$ for all $i\ge 1$ and all $m\ge 1$.
\item[{\rm (ii)}] The standard monomials, standard on $Y$ of degree $m$, form a basis
for $H^0(Y,\cl_{m\lam})$ for all $m\ge 1$.
\end{joliste}
\end{thm}

{\it Proof.}
We give only a sketch of the proof for $m=1$ (respectively $m=0$ in {\it i)}, the case
$m\ge 2$ can be proved in the same way as in the regular case. The proof of the appropriate
version of Proposition~\ref{smtonenonreg} for monomials of type $(s\rho,m\lam,q\rho)$ is the same.

The proof is (as in the regular case) by induction on the number of components of maximal
dimension  and on the dimension. The theorem holds obviously in case $Y$ is a point.
We assume now that the theorem holds for all pointed unions of Richardson varieties
of dimension smaller than $n$. Suppose $X_\kappa^\sigma$ is
of dimension $n$ and consider the projection $\phi: X_\kappa^\sigma \longrightarrow G/P_\lam$.

If $\delp_\lam X_\kappa^\sigma=\emptyset$,
then $X_\kappa^\sigma\subset p^{-1}(e_{\kbar})$ by Proposition~\ref{emptylambdaboundary}, and
the restriction of the line bundle to $X_\kappa^\sigma$ is trivial. The theorem holds obviously
in this case.

If $\delp_\lam X_\kappa^\sigma\not=\emptyset$, then the inclusion $\ci_{\delp_\lam X_\kappa^\sigma}
\hookrightarrow\co_{X_\kappa^\sigma}$, tensored by $\cl_\lam$, induces an exact sequence
$$
0\rightarrow \ci_{\delp_\lam X_\kappa^\sigma}(\lam) \rightarrow\co_{X_\kappa^\sigma}\otimes\cl_\lam
\rightarrow \co_{\delp_\lam X_\kappa^\sigma}\otimes\cl_\lam\rightarrow 0.
$$
By Theorem~\ref{idealfiltrationtwo}, the sheaf $\ci_{\delp_\lam X_\kappa^\sigma}(\lam)$
admits a filtration such that the associated graded is isomorphic to $\bigoplus \co_{X_{e(\pi)}^\sigma}
\otimes \chi_{-\pi(1)}$.

Since $\delp_\lam X_\kappa^\sigma$ is a pointed union of
Richardson varieties of smaller dimension, the vanishing theorem
for higher cohomology groups holds. The filtration in
Theorem~\ref{idealfiltrationtwo} and the vanishing in
Theorem~\ref{vanishandbasis}(i) shows that $\ci_{\delp_\lam
X_\kappa^\sigma}(\lam)$ has vanishing  higher cohomology groups,
so the long exact sequence shows that
$H^i(X_\kappa^\sigma,\cl_\lam)=0$ for all $i\ge 1$. Further, the
short exact sequence $$ 0\rightarrow
H^0(X_\kappa^\sigma,\ci_{\delp_\lam X_\kappa^\sigma}(\lam))
\rightarrow H^0(X_\kappa^\sigma,\cl_\lam) \rightarrow
H^0(\delp_\lam X_\kappa^\sigma,\cl_\lam)\rightarrow 0 $$ shows
$\dim H^0(X_\kappa^\sigma,\cl_\lam)=\dim H^0(\delp_\lam
X_\kappa^\sigma,\cl_\lam)+ \dim
H^0(X_\kappa^\sigma,\ci_{\delp_\lam X_\kappa^\sigma}(\lam))$.
Consider the right hand side of the equation. By induction, the
first term is equal to the  number of L-S paths of shape $\lam$
standard on $\delp_\lam X_\kappa^\sigma$, and the second is by the
filtration equal to the number of L-S paths of shape $\lam$
standard on $X_\kappa^\sigma$ and such that $i(\pi)=\kappa$. This
together is the number of L-S paths of shape $\lam$ standard on
$X_\kappa^\sigma$. Since we have already proved that the
restrictions of the corresponding path vectors remain linearly
independent (Proposition~\ref{smtonenonreg}), this proves that the
path vectors form a basis of $H^0(X_\kappa^\sigma,\cl_\lam)$.

Suppose now $Y$ is a pointed union of Richardson varieties. Let
$X_\kappa^\sigma$ be an irreducible component of maximal dimension
and denote $Y$ the union of the remaining irreducible components.
We have an exact sequence: $$ 0\rightarrow \co_{Y}\rightarrow
\co_{X_\kappa^\sigma}\oplus \co_{Y'}\rightarrow \co_{Y\cap
X_\kappa^\sigma}\rightarrow 0. $$ If we tensor the sequence with
$\cl_{\lam}$, then we know by induction on the dimension,
respectively the number of irreducible components of maximal
dimension that the higher cohomology groups of the second and the
last term vanish, and further, for the global sections the last
map is surjective. It follows that the higher cohomology groups
$H^i(Y,\cl_{\lam})$ vanish.

It remains to count the dimensions of the spaces of global
sections. Since one has standard monomial theory for the last two
terms by induction on the dimension, respectively the number of
irreducible components of maximal  dimension, one checks easily
that in fact $\dim H^0(Y,\cl_{\lam})$ is equal to the number of
standard monomials on $Y$. The linear independence
(Proposition~\ref{smtonenonreg}) then implies that the standard
monomials form a basis of $\dim H^0(Y,\cl_{\lam})$.
\endpf
\begin{rem}\rm
The result (i) in Theorem \ref{vanishandbasisnonreg} is also
proved in \cite{brionlak} (cf.  \cite{brionlak}, Proposition 1)
\end{rem}

\section{Standard monomials and equivariant
$K$-theory}\label{SMTequivariantKtheory}
In this section we will show that the notion of a ``standard monomial''
is directly related to geometry of Richardson varieties, and this can be best expressed
in terms of equivariant $K$-theory.

The Chow ring Chow$(G/B)$ has a $\bz$-basis consisting of
$[X_\tau]$, the class  in Chow$ (G/B)$ defined by the cycle
represented by the Schubert variety $X_\tau$, $\tau \in W$. Let
$\lambda \in \Lambda$ (=the weight lattice), $\cl_\lambda$ the
associated line bundle. The class $[\cl_{\lambda}]$ is defined by
a codimension one cycle. We have, the following expression for the
product $[X_\tau] \cdot [\cl_\lambda]$ in the Chow ring: $$
 [X_\tau] \cdot [\cl_\lambda]= \sum a_i [X_{\tau_i}]\leqno(1)
$$ where the summation runs over all Schubert divisors
$X_{\tau_i}$ in $X_\tau$; further, one has an explicit expression
for $a_i$ due to Chevalley (cf. \cite{che}). When $G = SL(n)$, (1)
is just the classical formula of Pieri.

Denote by $K(G/B)$ the Grothendieck ring of  (isomorphism classes
of) coherent sheaves on $G/B$.  Then one knows that the classes
$[{\mathcal{O}}_{X_\tau}]$ in $K(G/B)$ of the structure sheaves
${\mathcal{O}}_{X_\tau}$ of the Schubert varieties $X_\tau$ form a
$\bz$-basis of $K(G/B)$. Let $$
  [{\mathcal{O}}_{X_\tau}]\cdot [\cl_{\lambda}] = \sum_{\kappa \in W}
a_{\tau, \kappa}^{\lambda} [{\mathcal{O}}_{X_\kappa}]\leqno(2) $$
More generally, let $K_T (G/B)$ be the Grothendieck ring of
$T$-equivariant sheaves on $G/B$.  By \cite{KK} one knows again
that the classes $[{\mathcal{O}}_{X_\tau}]_T$ of the structure
sheaves of the Schubert varieties form a $\bz [\Lambda ]$-basis,
where ${\bz} [\Lambda ]$ is the group algebra of the weight
lattice $\Lambda$. Let $$
  [{\mathcal{O}}_{X_\tau}]_T\cdot [\cl_{\lambda}]_T = \sum_{\kappa \in W}
C_{\tau,\kappa}^{\lambda} [{\mathcal{O}}_{X_\kappa}]_T\leqno(3) $$
where $C_{\tau,\kappa}^{\lambda}$'s are formal sums of characters
of $T$ with integral coefficients. Assume now that $\lambda \in
\Lambda^+$. Then one knows that in (1) all $a_i \geq 0$
(cf.\cite{che}).  The integers $a_{\tau, \kappa}^\lambda$ in (2)
were determined  by Fulton and Lascoux (cf. \cite{FL}) for the
case $G = SL(n)$ ; they provide a formula for $a_{\tau,
\kappa}^{\lambda}$ (using the combinatorics of Grothendieck
polynomials). The general case was treated by Mathieu using
representation theory (cf. \cite{M}), who shows that
$C_{\tau,\kappa}^{\lambda}$ are effective i.e., formal sums of
characters with positive integral coefficients. For a L-S path
$\pi=(\tau_1> \cdots>\tau_r;a_1<\cdots <a_r)$ we set the initial
direction as $i(\pi)=\tau_1$ and the final direction
$e(\pi)=\tau_r$.

 We assume for simplicity that $\lam$ is {\it regular dominant}.
Using the path model theory, Pittie and Ram (cf. \cite{PR}) gave
an explicit determination of $C_{\tau,\kappa}^{\lambda}$'s: $$
[{\mathcal{O}}_{X_\tau}]_T\cdot [\cl_{\lambda}]_T =
\sum_{i(\pi)\le \tau} [{\mathcal{O}}_{X_{e(\pi)}}]_T\,
e^{\pi(1)}\leqno(4) $$ where the summation  runs over L-S paths
$\pi$ of shape $\lambda $.

An effective version of the relation (4) is proved in \cite{LiSe}
by constructing a filtration ${\cal{F}}:=\{{\cal{F}}^i\}$ for
${\cal{O}}_{X_\tau}\otimes \cl_\lambda$ by $B$-equivariant
${\cal{O}}_{G/B}$-modules such that each subquotient, as a
$B$-equivariant ${\cal{O}}_{G/B}$-sheaf, is isomorphic to
${\cal{O}}_{X_{e(\pi)}}\otimes \chi_{-\pi(1)}$, the structure
sheaf ${\cal{O}}_{X_{e(\pi)}}$, twisted by the character $-\pi(1)$
of $B$, for a suitable $\pi$ such that $i(\pi)\le \tau$.

Fixing a $\kappa\le \tau$, we obtain from (4) that the coefficient
of $[{\cal{O}}_{X_\kappa}]_T$ on the R.H.S. of (4) equals $$
\sum_{\{\pi|\tau\ge i(\pi), e(\pi)=\kappa\}}\ e^{\pi(1)}. $$ Using
the results of the preceding sections, we give this character a
representation theoretic interpretation. By tensoring the exact
sequence $$ 0\rightarrow \ci_{\delm X_\tau^\kappa} \rightarrow
{\mathcal{O}}_{X_\tau^\kappa}\rightarrow
{\mathcal{O}}_{\partial^-X_\tau^\kappa}\rightarrow 0 $$ by
$\cl_\lambda$, and writing the cohomology exact sequence, we get
the exact sequence $$ 0\rightarrow H^0(X_\tau^\kappa, \ci_{\delm
X_\tau^\kappa}(\lam) ) \rightarrow
H^0(X_\tau^\kappa,\cl_\lambda)\rightarrow
H^0(\partial^-(X_\tau^\kappa), \cl_\lambda) \rightarrow 0 $$ Hence
we obtain $\dim H^0(X_\tau^\kappa, \ci_{\delm X_\tau^\kappa}(\lam)
) = \# \{p_\pi | i(\pi)\le \tau,e(\pi)=\kappa\}$. Since these
vectors form a basis for the kernel of the second map, it follows
that $\{p_\pi | i(\pi)\le \tau,e(\pi)=\kappa\}$ is actually a
basis of $H^0(X_\tau^\kappa, \ci_{\delm X_\tau^\kappa}(\lam))$. In
the non-regular case, we have to work with the $\lambda$-boundary
$\partial_\lambda^-X_\tau^\kappa$ in the place of
$\partial^-X_\tau^\kappa$. Thus we obtain

\begin{thm}\label{k-smt}
With notations as above, let $\lambda$ be dominant. We have
$$C_{\tau,\kappa}^{\lambda}=\charc H^0(X_\tau^\kappa,
\cl_\lambda\otimes \ci_{\delm X_\tau^\kappa}), $$ $$
a_{\tau,\kappa}^{\lambda}=\#\{\pi \hbox{\rm{\ L--S path,\ shape\
}}\lambda \mid \tau\ge i(\pi),e(\pi)=\kappa \}. $$
\end{thm}

The above result is also proved in  \cite{brion} using a flat
family with generic fiber $\cong$  diag$(X_\tau)\subset
X_\tau\times X_\tau$, and the special fiber $\cong \cup_{x\le
\tau}X_x\times X_\tau^x$. We will give a different construction of
the flat family using SMT. Again the connection can be most
directly established in the language of equivariant $K$--theory.

Let $X_\tau^\kappa\subset G/Q$ be a Richardson variety. On the one
hand, consider the Richardson variety $X_\tau^\kappa$ diagonally
embedded in $Z=G/Q\times G/Q$, we write $\Delta X_\tau^\kappa$ for
this variety. On the other hand, consider the subvariety $$
Y=\bigcup_{\kappa\in W/W_Q\atop\kappa\le\sigma\le\tau}
X^\kappa_\sigma\times X^\sigma_\tau \,\subset X^\kappa_\tau\times
X^\kappa_\tau\,\subset Z. $$ We denote the corresponding structure
sheaves $\co_{\Delta X_\tau^\kappa}$, respectively $\co_Y$; we are
interested in describing their classes in the Grothendieck group
$K_T(Z)$ of $T$-equivariant coherent sheaves of $\co_Z$--modules
on $Z$.

\begin{thm}\label{k-smt1}
In $K_T(Z)$, the following equality of classes of $T$-equivariant coherent sheaves of $\co_{Z}$--modules holds:
$$
[\co_{\Delta \mbox{$\scriptstyle X_\tau^\kappa$}}]_T=[\co_{Y}]_T=\sum_{\sigma\in W/W_Q\atop\kappa\le\sigma\le\tau}
[\delp\ci_\sigma^\kappa\otimes\co_{X_\tau^\sigma}]_T=\sum_{\sigma\in W/W_Q\atop\kappa\le\sigma\le\tau}
[\co_{X_\sigma^\kappa}\otimes\delm\ci^\sigma_\tau]_T
$$
More precisely, there exist $T$-stable filtrations of $\co_{Y}$ as sheaf of $\co_Z$--modules, such that the
associated graded is the direct sum of the sheaves $\delp\ci_\sigma^\kappa\otimes\co_{X_\tau^\sigma}$,
respectively the direct sum of the sheaves $\co_{X_\sigma^\kappa}\otimes\delm\ci^\sigma_\tau$, where
$\kappa\le\sigma\le\tau$.

Further, there exists a flat family of subvarieties of $Z$ such
that the generic fibre is isomorphic to the diagonal Richardson
variety $\Delta X_\tau^\kappa$, and the special fibre is $Y$.
\end{thm}
\smallskip\noindent
{\it Proof.\/} We use again the correspondence between graded modules
over the homogeneous coordinate ring and coherent sheaves.
Fix a ($Q$-regular) dominant weight $\lam\in\Lambda_Q^{++}$ and let $B^\kappa_\tau(n\lam)$ be the
set of standard sequences of length $n$ of L--S paths of shape $\lam$, standard on $X^\kappa_\tau$.

Denote $\cR$ the homogeneous coordinate ring corresponding to the
embedding $Z\hookrightarrow\bp(H^0(Z,\cl_\lam\otimes\cl_\lam))$,
so that $$ \cR= \bigoplus_{n\ge 0}
H^0(Z,\cl_{n\lam}\otimes\cl_{n\lam}),$$ and let
$\cR_{\tau}^\kappa$ be the quotient: $$
\cR_{\tau}^\kappa=\bigoplus_{n\ge 0}
H^0(X_\tau^\kappa,\cl_{n\lam})\otimes
H^0(X_\tau^\kappa,\cl_{n\lam}) =\bigoplus_{n\ge 0}
H^0(X_\tau^\kappa\times
X_\tau^\kappa,\cl_{n\lam}\otimes\cl_{n\lam}). $$ Consider for
$n\ge 1$ the subspace (this is definitely not an ideal) $$
SM_\tau^\kappa(\lam)_n=\langle p_{\upi_1}\otimes p_{\upi_2}\mid
\upi_1, \upi_2\in B^\kappa_\tau(n\lam), \, p_{\upi_1}p_{\upi_2}
\hbox{\rm\ standard\ in}\, H^0(X_\tau^\kappa,\cl_{2n\lam})\rangle
$$ and set $SM_\tau^\kappa(\lam)_0=k$. Consider the graded vector
space: $$ SM_\tau^\kappa(\lam)=\bigoplus_{n\ge
0}SM_\tau^\kappa(\lam)_n. $$ We will present two ways to make this
graded vector space into a graded $\cR$-module.

First note that the subvariety $Y$ is a union of Richardson
varieties (for the group $G\times G$ with Borel subgroup $B\times
B$). So by Theorem~\ref{vanishandbasis}, the following restriction
map is $T$-equivariant, surjective, and maps the subspace
$SM_\tau^\kappa(\lam)$ isomorphically onto its image: $$
\cR_{\tau}^\kappa\rightarrow \cR_Y=\bigoplus_{n\ge 0}
H^0(Y,\cl_{n\lam}\otimes\cl_{n\lam}) $$ The kernel is in this case
also easy to describe: If the product $p_{\upi_1}p_{\upi_2}$ is
not standard, then $e(\upi_1)\not\ge i(\upi_2)$ and
$e(\upi_2)\not\ge i(\upi_1)$; hence either
$p_{\upi_1}\vert_{X_\sigma^\kappa}\equiv 0$ or
$p_{\upi_2}\vert_{X_\tau^\sigma}\equiv 0$. It follows that the
restriction of $p_{\upi_1}\otimes p_{\upi_2}$ on all irreducible
components $X_\sigma^\kappa\times X_\tau^\sigma$ vanishes. Since
the tensor products $p_{\upi_1}\otimes p_{\upi_2}$, $\upi_1,
\upi_2\in B_\tau^\kappa(n\lam)$, form a basis, it follows that the
kernel of the map is spanned by all $p_{\pi_1}\otimes p_{\pi_2}$
such that $\pi_1, \pi_2\in B_\tau^\kappa(n\lam)$ and
$p_{\upi_1}p_{\upi_2}$ is not standard in
$H^0(X_\tau^\kappa,\cl_{2n\lam})$.

Summarizing, by the $T$-equivariant graded vector space isomorphism of $SM_\tau^\kappa(\lam)\rightarrow\cR_Y$
we have endowed $SM_\tau^\kappa(\lam)$ in  a $T$-equivariant way with the structure of an $\cR$--module.

Next consider the product map $$ \cR_{\tau}^\kappa\rightarrow
\cR_{\mbox{$\scriptstyle \Delta X_\tau^\kappa$}} =\bigoplus_{n\ge
0} H^0(\Delta X_\tau^\kappa,\cl_{n\lam}\otimes\cl_{n\lam}) $$
Again by Theorem~\ref{vanishandbasis}, this map induces a
$T$-equivariant isomorphism of graded vector spaces
$SM_\tau^\kappa(\lam)\rightarrow\cR_{\mbox{$\scriptstyle \Delta
X_\tau^\kappa$}}$, so this induces a different structure as
$\cR$-module on $SM_\tau^\kappa(\lam)$.

Consider the $\bz[\Lambda]$-linear map
$\chi_\lam:K_T(Z)\longrightarrow \bz[\Lambda]$ defined on the
class of a $T$-equivariant coherent sheaf $\cf$ as follows: $$
\chi_\lam([\cf])=\sum_{i\ge 0} (-1)^i\charc
H^i(Z,\cf\otimes_{\co_{Z}}(\cl_\lam\otimes\cl_\lam)). $$ If two
elements of $K_T(Z)$ do not coincide, say $\sum
a_\cf[\cf]\not=\sum a_{\cf'}[\cf']$, then it is well known (see
for example~\cite{LiSe}) that there exists an $n\in\bn$ such that
$$ \chi_{n\lam}(\sum a_\cf[\cf])\not=\chi_{n\lam}(\sum
a_{\cf'}[\cf']). $$ Now by the vanishing of the higher cohomology
(Theorem~\ref{vanishandbasis}) and the $T$-equivariant graded
vector space isomorphisms above, we get: $$
\chi_{n\lam}(\co_Y)=\charc
SM_\tau^\kappa(\lam)_n=\chi_{n\lam}(\co_{\Delta X_\tau^\kappa}).
$$ and hence $[\co_Y]=[\co_{\Delta X_\tau^\kappa}]$ in $K_T(Z)$.

Next fix a numeration
$\kappa=\sigma_1,\sigma_2,\ldots,\tau=\sigma_q$ of the elements
between $\kappa$ and $\tau$ such that $\sigma_j>\sigma_i$ implies
$j>i$. Set $Y^j=\bigcup_{1\le i\le j} X^\kappa_{\sigma_i}\times
X^{\sigma_i}_\tau$, so we get a filtration $$
Y^1=X^\kappa_{\kappa}\times X^{\kappa}_\tau\subset
Y^2\subset\ldots\subset Y^q=Y. $$ Set $\cR_Y(j)=\bigoplus_{n\ge 0}
H^0(Y^j,\cl_{n\lam}\otimes \cl_{n\lam})$, this ring again has a
basis consisting of standard monomials. Set $K^{1}=\cR_Y$, and,
for $j=1,\ldots,q$, let $K^{j+1}\subset \cR_Y$ be the kernel of
the restriction map $\cR_Y\rightarrow \cR_Y(j)$. We get a
filtration $$ K^{q+1}=0\subset K^{q}\subset \ldots\subset
K^2\subset K^{1}=\cR_Y. $$ The kernels $K^j$ have again a basis by
standard monomials. In fact, recall the basis elements
$p_{\upi_1}\otimes p_{\upi_2}$ in $SM_\tau^\kappa(\lam)$, then
$K^j$  is spanned by those such that $i(\upi_1)=\sigma_i$ for some
$i\ge j$. Consider the map $$ K^j\rightarrow
H^0(X^\kappa_{\sigma_j}\times X^{\sigma_j}_\tau,\cl_{n\lam}\otimes
\cl_{n\lam}), \quad s\mapsto s\vert_{X^\kappa_{\sigma_j}\times
X^{\sigma_j}_\tau}. $$ It follows immediately that $K^{j+1}$ is
the kernel of this map, and, by the proof of Theorem~\ref{k-smt},
we get an isomorphism of $\cR$--modules: $$ K^j/K^{j+1}\simeq
\bigoplus_{n\ge 0} H^0(X^\kappa_{\sigma_j},\ci_{\delp
X^\kappa_{\sigma_j}}(n\lam)) \otimes
H^0(X^{\sigma_j}_\tau,\cl_{n\lam}) $$ Since all these maps are
$T$-equivariant, translated into the language of sheaves of
$\co_Z$--modules, this means that the filtration of $\cR_Y$
induces a $T$--stable filtration of $\co_Y$ as $\co_Z$--sheaf such
that the associated graded is a direct sum of the $T$--equivariant
$\co_Z$--sheaves $\ci_{\delp X^\kappa_{\sigma_j}}\otimes
\co_{X^{\sigma_j}_\tau}$.

The proof for the second filtration is similar and is left to the
reader. To describe the flat family note that the kernel of the
product map $$ \cR_{\tau}^\kappa\rightarrow
\cR_{\mbox{$\scriptstyle \Delta X_\tau^\kappa$}} =\bigoplus_{n\ge
0} H^0(\Delta X_\tau^\kappa,\cl_{n\lam}\otimes\cl_{n\lam}) $$ is
the ideal generated by the commutation relations $p_{\pi_1}\otimes
p_{\pi_2}-p_{\pi_2}\otimes p_{\pi_1}$ for all $\pi_1,\pi_2\in
B_\tau^\kappa(\lam)$, and if neither $p_{\pi_1}p_{\pi_2}$ nor
$p_{\pi_2}p_{\pi_1}$ is standard in
$H^0(X_\tau^\kappa,\cl_{2\lam})$, then we have the additional
relations $$ p_{\pi_1}\otimes p_{\pi_2} -\sum a_{\eta_1,\eta_2}
p_{\eta_1}\otimes p_{\eta_2} $$ where the $p_{\eta_1}p_{\eta_2}$
are standard and the coefficients $a_{\eta_1,\eta_2}\not=0$ only
if $(\eta_1,\eta_2)\succeq \pi_1\wedge\pi_2\succeq^r
(\eta_1,\eta_2)$. This follows easily from
Theorem~\ref{quadraticrelation}, using a monomial order as in
Proposition~\ref{homcoordinate}. Further, these elements form a
reduced Gr\"obner basis for the ideal. The existence of the flat
deformation follows now from standard Gr\"obner basis arguments.
\begin{rem}\rm
The existence of a filtration and a flat family (as in Theorem
\ref{k-smt1}) is also proved in \cite{brionlak}. The existence of
a filtration is proved in \cite{brion} for any Cohen-Macaulay
subvariety $X$ of a flag variety, in general position with respect
to opposite Schubert varieties (note that the assumptions on $X$
hold for Schubert varieties) . It is interesting to note that
while in this paper the existence of a filtration and a flat
family in Theorem \ref{k-smt1} is obtained as a consequence of
SMT, in \cite{brionlak} it is the starting point towards
developing a SMT !
\end{rem}


\goodbreak
\noindent
V.L.: Department of Mathematics, Northeastern University,
\par\noindent
Boston, MA 02115, USA
\vskip 8pt\noindent
P.L: Fachbereich Mathematik, Universit\"at Wuppertal,
Gau\ss-Stra\ss e 20
\par\noindent
42097 Wuppertal, Germany
\end{document}